\documentclass[10pt]{article}
\usepackage{latexsym,amsmath,amssymb,amscd}
\begin{document}

\makeatletter
\@addtoreset{figure}{section}
\def\thefigure{\thesection.\@arabic\c@figure}
\def\fps@figure{h,t}
\@addtoreset{table}{bsection}

\def\thetable{\thesection.\@arabic\c@table}
\def\fps@table{h, t}
\@addtoreset{equation}{section}
\def\theequation{\thesection.\arabic{equation}}
\makeatother

\newcommand{\bfi}{\bfseries\itshape}

\newtheorem{theorem}{Theorem}
\newtheorem{acknowledgment}[theorem]{Acknowledgment}
\newtheorem{algorithm}[theorem]{Algorithm}
\newtheorem{axiom}[theorem]{Axiom}
\newtheorem{case}[theorem]{Case}
\newtheorem{claim}[theorem]{Claim}
\newtheorem{conclusion}[theorem]{Conclusion}
\newtheorem{condition}[theorem]{Condition}
\newtheorem{conjecture}[theorem]{Conjecture}
\newtheorem{corollary}[theorem]{Corollary}
\newtheorem{criterion}[theorem]{Criterion}
\newtheorem{definition}[theorem]{Definition}
\newtheorem{example}[theorem]{Example}
\newtheorem{lemma}[theorem]{Lemma}
\newtheorem{notation}[theorem]{Notation}
\newtheorem{problem}[theorem]{Problem}
\newtheorem{proposition}[theorem]{Proposition}
\newtheorem{remark}[theorem]{Remark}
\numberwithin{theorem}{section}
\numberwithin{equation}{section}

%%% Todo
\newcommand{\todo}[1]{\vspace{5 mm}\par \noindent
\framebox{\begin{minipage}[c]{0.95 \textwidth}
\tt #1 \end{minipage}}\vspace{5 mm}\par}
%%%

\pagestyle{myheadings}
\markboth{\sl Belti\c t\u a and Ratiu: Symplectic leaves}
{\sl Belti\c t\u a and Ratiu: Symplectic leaves}

\makeatletter
\title{\textbf{Symplectic leaves in real Banach Lie-Poisson spaces}}
\author{Daniel Belti\c t\u a$^{1}$ and  Tudor S. Ratiu$^{2}$}
\addtocounter{footnote}{1}
\footnotetext{Institute of Mathematics ``Simion Stoilow'' of the Romanian
Academy,
RO-014700 Bucharest, Romania. \texttt{Daniel.Beltita@imar.ro}. }
\addtocounter{footnote}{1}
\footnotetext{Centre  Bernoulli,
\'Ecole Polytechnique F\'ed\'erale de Lausanne,  CH-1015 Lausanne,
Switzerland. \texttt{Tudor.Ratiu@epfl.ch}.}
\date{}
\makeatother
\maketitle

\begin{abstract}
We present several large classes of real Banach Lie-Poisson
spaces whose characteristic distributions are integrable,
the integral manifolds being symplectic leaves just as in
finite dimensions. We also investigate when these leaves
are embedded submanifolds or when they have K\"ahler
structures. Our results apply to the real Banach
Lie-Poisson spaces provided by the self-adjoint
parts of preduals of arbitrary $W^*$-algebras, as well
as of certain operator ideals.


\end{abstract}

\section{Introduction}

This paper studies some geometric properties of the recently
introduced  Banach Lie-Poisson spaces (see \cite{OR03}).
Every Banach Lie-Poisson space is the predual of
some Banach Lie algebra. Two classes of Banach Lie-Poisson
spaces will be investigated in this work:  preduals of $W^*$-algebras  and
preduals of certain operator ideals.

To explain the geometric questions addressed for these
two types of Banach Lie-Poisson spaces, recall that every
finite dimensional Poisson manifold has a characteristic generalized
distribution whose value at any point is the span of all
Hamiltonian vector fields evaluated at that point.  The
characteristic distribution is always integrable and each of
its leaves has two key features: it is an initial symplectic
submanifold of the Poisson manifold under consideration that is
at the same time a Poisson submanifold
(see e.g., \cite{W} or \cite{MR}).

If the
Poisson manifold is a Lie-Poisson space
${\mathfrak g}^*$, where ${\mathfrak g}$ is the Lie algebra of
some connected finite dimensional Lie group $G$,
it turns out that the integral manifolds of the characteristic
distribution  of ${\mathfrak g}^*$ are just the coadjoint orbits
of $G$  with the natural $G$-invariant orbit symplectic
structures (see e.g., \cite{W} or \cite{MR}).  If $G$ is compact, then the
coadjoint orbits are
$G$-homogeneous embedded K\"ahler submanifolds  of ${\mathfrak g}^*$ (see
e.g. \cite{GS}).

The goal of the present paper is to show that similar phenomena
occur in infinite dimensions for large classes of Banach
Lie-Poisson spaces.  The main results are described in
Corollaries ~\ref{corollary3.9}~and~\ref{corollary5.4} and in
Theorem~\ref{final}.
In the case of preduals of $W^*$-algebras,
weakly symplectic structures on integral manifolds of the
characteristic distribution  have been already constructed in
\cite{OR03} under a certain technical condition.  We shall prove
in Proposition~\ref{proposition3.7}  that this
condition is always satisfied,
for self-adjoint elements,
hence all the integral manifolds
of the corresponding characteristic distribution  are symplectic
leaves.  On the other hand, in the predual ${\mathfrak S}_1$
(trace class operators) of
${\mathcal B}({\mathcal H})$ (bounded operators)  for some complex Hilbert
space
${\mathcal H}$,  the question of which leaves are actually
embedded submanifolds of ${\mathfrak S}_1$  was answered in
\cite{Bo00} and \cite{Bo03}:  they are precisely the leaves containing
finite-rank operators.
It is noteworthy that a similar characterization of
the unitary orbits that are embedded submanifolds of
$\mathcal{B}(\mathcal{H})$
had been previously obtained in \cite{AS89},
cf. Theorem~\ref{theorem4.1} below.
(See \cite{AS91} for an extension of that characterization to
unitary orbits in arbitrary $C^*$-algebras.)
We shall prove a similar result in
the more general setting of operator ideals  (see
Theorem~\ref{final}).  Moreover, we will show that all
these embedded submanifolds are actually weakly K\"ahler
homogeneous spaces,  thus recovering what happens in finite
dimensions for the coadjoint orbits  of the compact group
${\rm U}(n)$. This circle of ideas is naturally related to
the general question of prequantization of
infinite dimensional manifolds carrying a closed two form 
and the problem of finding Banach Lie groups acting naturally
on the  relevant associated bundles; see \cite{Ne02, Ne03, Ne03a}
for progress in this direction.

\section{Symplectic leaves in preduals of $W^*$-algebras}

Throughout the paper, by $C^*$-algebra (respectively $W^*$-algebra)
we actually mean \textit{unital} $C^*$-algebra
(respectively \textit{unital} $W^*$-algebra).

\begin{definition}\label{notation3.1}
%\normalfont
For every $C^*$-algebra $M$, let
$$\mathcal{P}_M:=\{p\in M\mid p^2=p^*=p\}$$
be the set of all orthogonal projections in $M$.
We denote by ${\operatorname{U}}_M$ the Banach Lie group of
all unitary elements of $M$.
Every $u\in {\operatorname{U}}_M$ defines an isometric $*$-isomorphism
$$\operatorname{Ad}(u)\colon M\to M,\quad a\mapsto uau^*.$$
We also denote by ${\mathfrak u}_M$ the Lie algebra of 
${\operatorname{U}}_M$, that is,
$${\mathfrak u}_M:=\{a\in M\mid a^*=-a\}.$$
\end{definition}

Throughout the paper, if $M$ is a $W^*$-algebra then
$M_* $ denotes the predual of  $M $ and
$M^* $ the dual of $M $.
An element $\varphi \in M_* $ is said
to be \textit{faithful\/} if
$\varphi(a^*a)>0$ whenever $0\ne a\in M$.
This condition is equivalent to the fact that the support of $\varphi$
equals $\mathbf{1}$
(see Remark~\ref{remark3.4} below).

Recall that a smooth map $f\colon V\to W $ between the
Banach manifolds $V $ and $W$ is said to be a
\textit{weak immersion\/}
if its tangent map $T_v f \colon T_v V\to T_{f(v)} W $
at any point $v \in V $ is an injective linear
bounded map.
Note that no assumption about  the closedness of
the range and its splitting properties are made.

\begin{theorem}\label{theorem3.2}
{{\rm (\cite{AV02})}}
Let $M$ be a $W^*$-algebra and $\varphi\in M_*$ faithful.
Consider the centralizer of $\varphi$, that is, the
sub-$W^*$-algebra
$$M^\varphi=\{a\in M\mid (\forall b\in M)\quad
\varphi(ab)=\varphi(ba)\},$$
and the unitary orbit of $\varphi$,
$$\mathcal{U}_\varphi=\{\varphi\circ\operatorname{Ad}(u)\mid u\in 
{\operatorname{U}}_M\}\simeq
{\operatorname{U}}_M/{\operatorname{U}}_{M^\varphi},$$
where ${\operatorname{U}}_{M^ \varphi}:= \{a \in {\operatorname{U}}_M \mid 
(\forall b\in M)\quad
\varphi(ab)=\varphi(ba)\}$ is
the unitary group of the centralizer algebra
$M^\varphi$, that is, the unitary elements of $M^ \varphi$.

Then the following assertions hold:

\begin{itemize}
\item[{\rm(i)}] $\mathcal{U}_\varphi\subseteq M_*$.
\item[{\rm(ii)}] The unitary group ${\operatorname{U}}_{M^\varphi}$ of the 
centralizer 
algebra $M^\varphi$
is a Lie subgroup of ${\operatorname{U}}_M$.
\item[{\rm(iii)}] The unitary orbit $\mathcal{U}_\varphi$
has a natural structure of weakly immersed submanifold of $M_*$
and ${\operatorname{U}}_M $ acts on it smoothly on the left via $(u, \psi) 
\in {\operatorname{U}}_M
\times
\mathcal{U}_\varphi \mapsto  \psi \circ \operatorname{Ad}(u
^{-1}) \in \mathcal{U}_\varphi$.
\item[{\rm(iv)}] The smooth manifold $\mathcal{U}_\varphi$ is simply
connected.
\end{itemize}
\end{theorem}

\noindent\textbf{Proof.\ \ }
{\rm(i)} This is obvious.

{\rm(ii)} Note that ${\operatorname{U}}_{M^\varphi}$ is an algebraic 
subgroup of
${\operatorname{U}}_M$
in the following sense
(see Definition~8.9 in \cite{Be02}):
$${\operatorname{U}}_{M^\varphi}=\{a\in {\operatorname{U}}_M\mid (\forall 
p\in{\mathcal P})\quad
p(a,a^{-1})=0\},$$
where ${\mathcal P}$ is a set of continuous polynomial functions on
$M\times M$.
In fact, we may take ${\mathcal P}=\{p_b\}_{b\in M}$,
where
$$p_b\colon M\times M\to{\mathbb C},\quad
p_b(x,y)=\varphi(xb)-\varphi(bx)$$
whenever $b\in M$; note that the polynomial $p_b$ depends only
on $x $, but we think of it as a function of $(x, y)$.
It is clear that each $p_b$ is a continuous linear functional on $M\times
M$,
and thus a polynomial of degree $\le1$.

Now the fact that ${\operatorname{U}}_{M^\varphi}$
is a Lie group with the topology inherited from ${\operatorname{U}}_M$  
follows by
the main result of
\cite{HK77};
see Theorem~8.12 in \cite{Be02} for the
precise statement in this regard.
Furthermore, to prove
that
${\operatorname{U}}_{M^\varphi}$ is actually a Lie subgroup of 
${\operatorname{U}}_M$,  we still
have to show that the Lie algebra ${\mathfrak u}_{M^\varphi}$
is a split subspace of ${\mathfrak u}_M$.  The latter fact is a
consequence of the fact that, since $\varphi$ is a normal
faithful positive form  on
$M$, there exists a conditional expectation
$E$ of $M$ onto $M^\varphi$.
We recall from Lemma~8.14.6 in \cite{P} that $M^\varphi$ equals 
the fixed-point algebra of the modular group of automorphisms
of $M$ associated with $\varphi$.
Thus the main theorem of \cite{Ta72} implies that there exists 
a conditional expectation $E$ from $M$ onto $M^\varphi$
satisfying $\varphi\circ E=\varphi$.
(Alternatively, the existence of $E$ follows by Remark~2.1 in 
\cite{AV99}.) 

{\rm(iii)} The unitary orbit
$\mathcal{U}_\varphi=\{\varphi\circ\operatorname{Ad}(u)\mid u\in
{\operatorname{U}}_M\}$ through
$\varphi \in M_*$ is in bijective correspondence with
${\operatorname{U}}_M/({\operatorname{U}}_{M})_\varphi$, where
\[
({\operatorname{U}}_M)_\varphi : = \{ u \in {\operatorname{U}}_M \mid
\varphi\circ\operatorname{Ad}(u) =
\varphi\}
\]
is the isotropy subgroup of $\varphi$ under the dual of the
action $\operatorname{Ad}$, where
$\operatorname{Ad}(u) b : = ubu ^{-1}$ for any $b \in M $. It
is easily verified that
$$({\operatorname{U}}_M)_\varphi = {\operatorname{U}}_{M^ \varphi} .$$
By
(ii), ${\operatorname{U}}_{M^ \varphi} $ is a Lie subgroup of 
${\operatorname{U}}_M $
and thus
the set ${\operatorname{U}}_M/{\operatorname{U}}_{M^ \varphi} $
has a unique smooth manifold
structure making the canonical projection ${\operatorname{U}}_M 
\rightarrow
{\operatorname{U}}_M/{\operatorname{U}}_{M^ \varphi}$ a surjective 
submersion;
the underlying
manifold topology of ${\operatorname{U}}_M/{\operatorname{U}}_{M^ \varphi} 
$
is the quotient
topology and ${\operatorname{U}}_M $ acts smoothly on the left on 
${\operatorname{U}}_M/{\operatorname{U}}_{M^
\varphi} $ by $(u, [v]) \in {\operatorname{U}}_M \times
{\operatorname{U}}_M/{\operatorname{U}}_{M^ \varphi}
\mapsto u\cdot [v]: = [uv]$, where $[v] = v{\operatorname{U}}_{M^ 
\varphi}$
(see
Bourbaki \cite{Bo3}, Chapter III, \S 1.6, Proposition 11).
Endow
the orbit $\mathcal{U}_\varphi$ with the manifold structure
making  this equivariant bijection into a diffeomorphism. It is
then easily checked that the inclusion of $\mathcal{U}_\varphi$
into $M_* $ is a weak immersion.

{\rm(iv)} See Theorem~2.9 in \cite{AV02}.
\quad $\blacksquare$

\begin{remark}\label{remark3.3a}
{\rm(cf.~Remark~A.2.2 in \cite{JS97})}
\normalfont
There always exist faithful elements in $M_*$ provided the predual $M_*$
of the $W^*$-algebra $M$ is separable.
\end{remark}

\begin{remark}\label{remark3.3b}
{\rm(cf.~Proposition~5.1 in \cite{AV99})}
\normalfont
In the setting of Theorem~\ref{theorem3.2}, assume that $M=\mathcal{ 
B}(\mathcal{ H})$
for some complex infinite dimensional Hilbert space $\mathcal{
H}$.

Then, for any faithful state $\varphi\in M_*$ the orbit
$\mathcal{U}_\varphi$ is not locally closed in
$M_*$. Thus, if  $M=\mathcal{ B}(\mathcal{ H})$,
the weakly immersed submanifolds occurring in
Theorem~\ref{theorem3.2}
are never embedded submanifolds of $M_*$.
\end{remark}

\begin{remark}\label{remark3.4}
{\rm(cf.~Section~5.15 in \cite{SZ79})}
\normalfont
Let $M$ be a $W^*$-algebra and $0\le\varphi\in M_*$.
Define the {\it support} of $\varphi$ by
$$p:=\mathbf{s}(\varphi):=1-\sup\{q\in\mathcal{ P}_M\mid
\varphi(q)=0\}\in\mathcal{ P}_M.$$
The support of $\varphi$ has the following  properties:

\begin{itemize}
\item[{\rm(i)}] $(\forall x\in M)\quad \varphi(x)
=\varphi(xp)=\varphi(px)
=\varphi(pxp)$.
\item[{\rm(ii)}] If $0\le x\in M$ and $\varphi(x)=0$ then
$pxp=0$.  In particular,
$\varphi|_{pMp}\in(pMp)_*$
is faithful on the $W^*$-algebra $pMp$.
\end{itemize}

For later reference we also note that we have
\begin{equation}\label{*}
(\forall u\in {\operatorname{U}}_M)\quad 
\mathbf{s}(\operatorname{Ad}(u)^*\varphi)
=u^{-1}\mathbf{s}(\varphi)u,
\end{equation}
since for each $q\in\mathcal{ P}_M$ the condition $\varphi(uqu^{-1})=0$
is equivalent to $uqu^{-1}\le1-\mathbf{s}(\varphi)$, hence to
$q\le1-u^{-1}\mathbf{s}(\varphi)u$.
\end{remark}

\begin{remark}\label{remark3.5}
{\rm(cf.~Section~5.17 in \cite{SZ79})}
\normalfont
Let $M$ be a $W^*$-algebra and $\varphi\in M_*$ such that
$\varphi=\varphi^*$,
in the sense that $\varphi(x^*)=\overline{\varphi(x)}$ for all $x\in M$.
Then there exist $\varphi_1,\varphi_2\in M_*$ uniquely determined by the 
conditions:
\begin{itemize}
\item[{\rm(i)}] $\varphi=\varphi_1-\varphi_2$,

\item[{\rm(ii)}] $\varphi_1\ge0$ and $\varphi_2\ge0$, and

\item[{\rm(iii)}] $\mathbf{s}(\varphi_1)\mathbf{s}(\varphi_2)=0$.

\end{itemize}

\end{remark}

\begin{lemma}\label{lemma3.6}
Let $M$ be a $W^*$-algebra, $0\le\varphi\in M_*$, 
$p:=\mathbf{s}(\varphi)$,
$\varphi_p:=\varphi|_{pMp}\in(pMp)_*$,
and denote, as before,
$${\operatorname{U}}_{M^\varphi}=\{u\in {\operatorname{U}}_M\mid 
\operatorname{Ad}(u)^*\varphi=\varphi\}.$$
Then
$$\begin{aligned}
{\operatorname{U}}_{M^\varphi}&=\{u\in {\operatorname{U}}_M\mid pu
=up,\ pup\in {\operatorname{U}}_{(pMp)^{\varphi_p}}\} \\
&= \left\{
\left(
\begin{array}{cc}
u_1 & 0\\
0 & u_2
\end{array}
\right)
\in {\operatorname{U}}_p
\;\Big|\; u_1\in {\operatorname{U}}_{(pMp)^{\varphi_p}},\, u_2\in
{\operatorname{U}}_{(1-p)M(1-p)}\right\},
\end{aligned}$$
where
$${\operatorname{U}}_p:=\{u\in {\operatorname{U}}_M\mid pu=up\}
= \left\{\left(
\begin{array}{cc}
u_1 & 0\\
0 & u_2
\end{array}
\right)
\;\Big|\;  u_1\in {\operatorname{U}}_{pMp}, u_2\in
{\operatorname{U}}_{(1-p)M(1-p)}\right\}$$
and the $2\times2$ matrix is written with respect to the
orthogonal decomposition $1=p+(1-p)$.
\end{lemma}

\noindent\textbf{Proof.\ \ }
Since
$$\{u\in {\operatorname{U}}_M\mid pu=up\}
= \left\{\left(
\begin{array}{cc}
u_1 & 0\\
0 & u_2
\end{array}
\right)
\;\Big|\;  u_1\in {\operatorname{U}}_{pMp}, u_2\in
{\operatorname{U}}_{(1-p)M(1-p)}\right\}$$
it follows that ${\operatorname{U}}_p$ is a Lie subgroup of
${\operatorname{U}}_M$.  For all $u\in {\operatorname{U}}_{M^\varphi}$ we 
have
$u^{-1}\mathbf{s}(\varphi)u=\mathbf{s}(\varphi)$  by formula~\eqref{*} in
Remark~\ref{remark3.4}.  Thus, since $p=\mathbf{s}(\varphi)$, we get
$${\operatorname{U}}_{M^\varphi}\subseteq {\operatorname{U}}_p.$$
We now come back to the proof of the desired conclusion.
For any $u\in {\operatorname{U}}_M$ we have
$$(\forall x\in M)\quad \varphi(uxu^{-1})=\varphi(x)
\iff (\forall x\in M)\quad \varphi(puxu^{-1}p)=\varphi(pxp)$$
by Remark~\ref{remark3.4}(i).
Hence for $u\in {\operatorname{U}}_p$ (that is, $up=pu$) we have
$$(\forall x\in M)\quad \varphi(uxu^{-1})=\varphi(x)
\iff (\forall x\in M)\quad \varphi((pup)(pxp)(pu^{-1}p))=\varphi(pxp).$$
Next note that, since $up=pu$, it follows that $pu^{-1}p$
is just the inverse of $u_1:=pup\in {\operatorname{U}}_{pMp}$.
Thus the above equivalence shows that,
for $u=\left(
\begin{array}{cc}
u_1 & 0\\
0 & u_2
\end{array}
\right)
\in {\operatorname{U}}_p$ as above,
we have
$$u\in {\operatorname{U}}_{M^\varphi}\iff u_1\in
{\operatorname{U}}_{(pMp)^{\varphi_p}},$$
and the desired conclusion is proved.
\quad $\blacksquare$

\begin{proposition}\label{proposition3.7}
Let $M$ be a $W^*$-algebra, $\varphi\in M_*$ such that
$\varphi=\varphi^*$, and
$${\operatorname{U}}_{M^\varphi}=\{u\in {\operatorname{U}}_M\mid 
\operatorname{Ad}(u)^*\varphi=\varphi\}.$$
Then ${\operatorname{U}}_{M^\varphi}$
is a Lie subgroup of ${\operatorname{U}}_M$.
\end{proposition}

\noindent\textbf{Proof.\ \ }
Let $\varphi=\varphi_1-\varphi_2$ as in Remark~\ref{remark3.5},
and denote $p_1=\mathbf{s}(\varphi_1)$, $p_2=\mathbf{s}(\varphi_2)$,
so that $p_1p_2=p_2p_1=0$.
We will prove that
$${\operatorname{U}}_{M^\varphi}={\operatorname{U}}_{M^{\varphi_1}}\cap 
{\operatorname{U}}_{M^{\varphi_2}}.$$
The inclusion $\supseteq$ is obvious.
Now let $u\in {\operatorname{U}}_{M^\varphi}$.
Then
$$\varphi=\operatorname{Ad}(u)^*\varphi=
\operatorname{Ad}(u)^*\varphi_1-\operatorname{Ad}(u)^*\varphi_2.$$
Moreover, it is clear that $\operatorname{Ad}(u)^*\varphi_j\ge0$ and
$\mathbf{s}(\operatorname{Ad}(u)^*\varphi_j)=u^{-1}p_ju$ (by \eqref{*} in 
Remark~\ref{remark3.4}) for $j=1,2$,
hence $\mathbf{s}(\operatorname{Ad}(u)^*\varphi_1)\mathbf{s}
(\operatorname{Ad}(u)^*\varphi_2)=0$.
It then follows from the uniqueness assertion in Remark~\ref{remark3.5} 
that
$\operatorname{Ad}(u)^*\varphi_j=\varphi_j$ for $j=1,2$, hence
$u\in {\operatorname{U}}_{M^{\varphi_1}}\cap 
{\operatorname{U}}_{M^{\varphi_2}}$
as desired.

Next denote $p_3=1-p_1-p_2$, so that $p_j\in \mathcal{ P}_M$ and 
$p_ip_j=0$
for $1\le i,j\le 3$, and $p_1+p_2+p_3=1$.
According to Lemma~\ref{lemma3.6} we have
${\operatorname{U}}_{M^{\varphi_{p_j}}}\subseteq\{u\in
{\operatorname{U}}_M\mid up_j=p_ju\}$
for $j=1,2$,
hence
$$\aligned
{\operatorname{U}}_{M^{\varphi_1}}\cap 
{\operatorname{U}}_{M^{\varphi_2}}&\subseteq\{u\in 
{\operatorname{U}}_M\mid
up_j=p_ju
\text{ for }j=1,2,3\} \\
&\simeq
\left\{\left(
\begin{array}{ccc}
u_1 & 0 & 0\\
0 & u_3 & 0 \\
 0 & 0  &u_2
\end{array}
\right)
\;\Big|\; u_j\in {\operatorname{U}}_{p_jMp_j}\text{ for }j=1,2,3\right\}.
\endaligned$$
Lemma~\ref{lemma3.6} actually shows that
$$\aligned
{\operatorname{U}}_{M^{\varphi_{p_1}}}
&=\{u\in {\operatorname{U}}_M\mid p_1u=up_1,\ p_1up_1\in 
{\operatorname{U}}_{(p_1Mp_1)^{\varphi_{p_1}}}\}\\
&\simeq {\operatorname{U}}_{M^{\varphi_{p_1}}}\times 
{\operatorname{U}}_{(1-p_1)M(1-p_2)}\\
&={\operatorname{U}}_{M^{\varphi_{p_1}}}\times 
{\operatorname{U}}_{(p_3+p_2)M(p_3+p_2)}
\endaligned$$
and similarly
$$\aligned
{\operatorname{U}}_{M^{\varphi_{p_2}}}
&=\{u\in {\operatorname{U}}_M\mid p_2u=up_2,\ p_2up_2\in
{\operatorname{U}}_{(p_2Mp_2)^{\varphi_{p_2}}}\}\\
&\simeq {\operatorname{U}}_{(p_1+p_3)M(p_1+p_3)}\times 
{\operatorname{U}}_{M^{\varphi_{p_2}}}.
\endaligned$$
Hence
$$\aligned
{\operatorname{U}}_{M^{\varphi_1}}\cap {\operatorname{U}}_{M^{\varphi_2}}
&=\{u\in {\operatorname{U}}_M\mid p_ju=up_j,\ p_jup_j\in 
{\operatorname{U}}_{(p_jMp_j)^{\varphi_{p_j}}}\text{ for }j=1,2\}\\
&\simeq {\operatorname{U}}_{M^{\varphi_{p_1}}} \times 
{\operatorname{U}}_{p_3Mp_3} \times 
{\operatorname{U}}_{M^{\varphi_{p_2}}} \\
&\simeq
\left\{\left(
\begin{array}{ccc}
u_1 & 0 & 0\\
0 & u_3 & 0 \\
 0 & 0  &u_2
\end{array}
\right) \;\Big|\; u_3\in {\operatorname{U}}_{p_3Mp_2},\ u_j\in
{\operatorname{U}}_{M^{\varphi_{p_j}}}
                    \text{ for }j=1,2\right\}.
\endaligned$$
Now ${\operatorname{U}}_{M^{\varphi_{p_j}}}$ is a Lie subgroup
of ${\operatorname{U}}_{p_jMp_j}$ by
Theorem~\ref{theorem3.2}(ii) since $\varphi_{p_j}=\varphi|_{p_jMp_j}$ is 
faithful for $j=1,2$
by Remark~\ref{remark3.4}(ii).
Hence the above isomorphism shows that 
${\operatorname{U}}_{M^{\varphi_1}}\cap
{\operatorname{U}}_{M^{\varphi_2}}$
is a Lie subgroup of
${\operatorname{U}}_{p_1Mp_1}\times {\operatorname{U}}_{p_3Mp_3}\times 
{\operatorname{U}}_{p_2Mp_2}$.
But the latter group is isomorphic to
$$\left\{\left(
\begin{array}{ccc}
u_1 & 0 & 0\\
0 & u_3 & 0 \\
 0 & 0  &u_2
\end{array}
\right) \;\Big|\; u_j\in {\operatorname{U}}_{p_jMp_j}\text{ for 
}j=1,2,3\right\},$$
which is a Lie subgroup of ${\operatorname{U}}_M$,
hence 
${\operatorname{U}}_{M^\varphi}={\operatorname{U}}_{M^{\varphi_1}}\cap
{\operatorname{U}}_{M^{\varphi_2}}$
is in turn a Lie subgroup of ${\operatorname{U}}_M$.
\quad $\blacksquare$

\begin{corollary}\label{corollary3.8}
For every $W^*$-algebra $M$ and $\varphi=\varphi^*\in M_*$,
the coadjoint orbit of the Lie group ${\operatorname{U}}_M$ through
$\varphi \in({\mathfrak u}_M)_*\subseteq({\mathfrak u}_M)^*$
has the structure of a ${\operatorname{U}}_M$-homogeneous weakly
symplectic
manifold which is weakly immersed into
$({\mathfrak u}_M)_*$.
\end{corollary}

\noindent\textbf{Proof.\ \ }
Just use Proposition~\ref{proposition3.7} along with Theorem~7.3 in 
\cite{OR03}.
\quad $\blacksquare$

\begin{corollary}\label{corollary3.9}
Let $M$ be an arbitrary $W^*$-algebra and
consider the corresponding real Banach Lie-Poisson space
$M_*^{\text{sa}}=\{\varphi\in M_*\mid \varphi=\varphi^*\}$.
Then the characteristic distribution of $M_*^{\text{sa}}$ is integrable 
and
all its maximal integral manifolds are symplectic leaves.
\end{corollary}

\noindent\textbf{Proof.\ \ }
Use Corollary~\ref{corollary3.8} along with Theorem~7.4 in
\cite{OR03} and note that all the coadjoint orbits referred to
in Corollary~\ref{corollary3.8}
 are connected since the unitary group of every $W^*$-algebra is
connected.
\quad $\blacksquare$
\medskip

\begin{remark}\label{strong}
\normalfont
It is noteworthy that the weakly symplectic manifolds given by
Corollary~\ref{corollary3.8} are sometimes strongly symplectic.
For instance, this is the case of the coadjoint orbits of rank-one projections 
if we assume  that $M=\mathcal{B}(\mathcal{H})$ for some complex 
Hilbert space ${\mathcal H}$ with the scalar product $(\cdot\mid\cdot)$. 

In fact,  for any $x\in{\mathcal H}$ with $\|x\|=1$ 
denote by $p_x=(\cdot\mid x)x$ the orthogonal projection of ${\mathcal H}$ 
onto the one-dimensional subspace ${\mathbb C}x$. 
Then $p_x\in M_*^{\text{sa}}$ and $up_xu^*=p_{ux}$ 
for all unit vectors $x\in{\mathcal H}$ and all $u\in{\operatorname{U}}_M$. 
Thus, denoting by $S_{\mathcal H}$ the unit sphere of ${\mathcal H}$
(that is, the set of all unit vectors in ${\mathcal H}$) 
and by 
${\mathbb P}({\mathcal H}):=S_{\mathcal H}/{\mathbb T}$ 
the projective space of ${\mathcal H}$, 
it follows that the mapping
$$S_{\mathcal H}\to M_*^{\text{sa}},\quad x\mapsto p_x,$$
induces a ${\operatorname{U}}_M$-equivariant diffeomorphism of 
${\mathbb P}({\mathcal H})$ onto the coadjoint orbit 
${\mathcal P}_1:=\{p_x\mid x\in S_{\mathcal H}\}$.  
It is well known that the projective space ${\mathbb P}({\mathcal H})$ is 
a strongly symplectic manifold (it is locally symplectomorphic
to ${\mathcal H}/{\mathbb C}x_0$ with the symplectic form defined by 
the double of the imaginary part of the quotient scalar product, 
for an arbitrary unit vector $x_0\in{\mathcal H}$), 
hence our claim that ${\mathcal P}_1$ is strongly symplectic will follow as soon as we
show that the aforementioned diffeomorphism 
${\mathcal P}_1\simeq{\mathbb P}({\mathcal H})$ 
is actually a symplectomorphism. 
To this end, fix a unit vector $x_0\in{\mathcal H}$. 
The symplectic structure of the coadjoint orbit 
${\mathcal P}_1$ through $p_{x_0}\in M_*^{\text{sa}}$ 
is defined by the skew symmetric bilinear form 
$$\omega_{x_0}\colon{\mathfrak u}_M\times{\mathfrak u}_M\to{\mathbb R},\quad 
\omega_{x_0}(a_1,a_2)=i\text{Tr}\,(p_{x_0}[a_1,a_2])$$
(see formula~(7.5) in \cite{OR03}).
Since the elements of ${\mathfrak u}_M$ are skew-symmetric,
it follows that for all $a_1,a_2\in{\mathfrak u}_M$ we have
$$\omega_{x_0}(a_1,a_2)
=-i\cdot\text{Tr}\,((\cdot\mid[a_1,a_2]x_0)x_0)
=i([a_1,a_2]x_0\mid x_0)
=2\text{Im}\,(a_1x_0\mid a_2x_0).
$$
On the other hand, if we consider
$\operatorname{U}_M^{x_0}=\{u\in\operatorname{U}\mid ux_0\in{\mathbb C}x_0\}$,
which is the isotropy group of ${\mathbb C}x_0\in{\mathbb P}({\mathcal H})$,
then we have a $\operatorname{U}_M$-equivariant diffeomorphism
$\operatorname{U}_M/\operatorname{U}_M^{x_0}\simeq{\mathbb P}({\mathcal H})$,
and the $\operatorname{U}_M$-invariant symplectic form of ${\mathbb P}({\mathcal H})$
will be defined by the skew-symmetric bilinear form
$$\omega'_{x_0}\colon{\mathfrak u}_M\times{\mathfrak u}_M\to{\mathbb R},\quad
\omega_{x_0}(a_1,a_2)=2\text{Im}\,(a_1x_0\mid a_2x_0).$$
The above computation shows that $\omega_{x_0}=\omega'_{x_0}$,
and this concludes the proof of the fact that the $\operatorname{U}_M$-equivariant
diffeomorphism ${\mathcal P}_1\simeq{\mathbb P}({\mathcal H})$
is a symplectomorphism, whence the coadjoint orbit ${\mathcal P}_1$ is
strongly symplectic.
\end{remark}

In the same special case when $M=\mathcal{B}(\mathcal{H})$ for some complex 
Hilbert space 
${\mathcal H}$ (that is, when $M$ is a factor of type~$\mbox{\rm I}$),
the result of Corollary~\ref{corollary3.8} also follows by 
Corollary~7.7 in \cite{OR03} along with Lemma~4.1 in \cite{AV99}.
However, we conclude this section by showing that, in general, the 
previous 
Corollary~\ref{corollary3.8}
applies to coadjoint orbits that do not fall under the hypotheses of 
Corollary~7.7 in \cite{OR03}.
To this end, we prove the following fact.

\begin{proposition}\label{new}
Let $M$ be a von Neumann algebra on a complex Hilbert space ${\mathcal
H}$.
Assume that $M$ is a $\mbox{\rm II}_1$ factor with the faithful normal 
trace
state $\tau$
and that there is a positive invertible element $h\in M$ with the spectral
measure $E_h(\cdot)$
such that for some $v\in{\mathcal H}\setminus\{0\}$ the localized measure
$\|E_h(\cdot)v\|^2$ has no atoms.

Now define
$$\varphi\colon M\to{\mathbb C},\quad \varphi(x)=\tau(hx).$$
Then $\varphi\in M_*$ is a faithful functional and
there exists no family $\{e_i\}_{i\in I}$ of mutually orthogonal
self-adjoint projections in $M$ satisfying $\sum\limits_{i\in 
I}e_i=\mathbf{1}$
and $M^{\varphi}=\Bigl\{\sum\limits_{i\in I}e_ixe_i\mid x\in M\Bigr\}$.
\end{proposition}

\noindent\textbf{Proof.\ \ }
It is clear that $\varphi\in M_*$.
Next, for every $x\in M$ we have
$$\varphi(x^*x)=\tau(hx^*x)=\tau(xhx^*)=\tau((h^{1/2}x^*)^*(h^{1/2}x^*)).$$
Since $\tau$ is faithful and $h$ is invertible,
it then easily follows that $\varphi$ is faithful.

Now, to prove the property stated for $M^\varphi$, we first check that
$$M^\varphi=\{a\in M\mid ah=ha\}.$$
In fact, $a\in M^\varphi$ if and only if $\varphi(ax)=\varphi(xa)$ for all 
$x\in M$,
that is, $\tau(hax)=\tau(hxa)$ for $x\in M$.
Since $\tau$ is a trace, the latter property is further equivalent to
$\tau(hax)=\tau(ahx)$ for all $x\in M$, hence to $ha=ah$.

Next let us assume that there exists a family $\{e_i\}_{i\in I}$ of
self-adjoint projections in $M$ satisfying $e_ie_j=0$ whenever $i\ne j$,
$\sum\limits_{i\in I}e_i=\mathbf{1}$ 
and $M^{\varphi}=\Bigl\{\sum\limits_{i\in I}e_ixe_i\mid x\in M\Bigr\}$.
According to the previous characterization of $M^\varphi$
we have that $h$ belongs to the center of $M^\varphi$.
Then it follows at once that $e_ihe_i$ ($=e_ih=he_i$)
belongs to the center of $e_iMe_i$ for each $i\in I$.
On the other hand, since $M$ is a factor
(i.e., its center reduces to the scalar multiples of the unit element)
it follows by Corollary~3.15 in \cite{SZ79} that $e_iMe_i$ is in turn a 
factor,
hence there exists $\lambda_i\in{\mathbb C}$ such that 
$e_ih=he_i=\lambda_ie_i$,
for arbitrary $i\in I$.

We now show that this fact contradicts the spectral assumption on $h$.
In fact, since the measure $\|E_h(\cdot)v\|^2$ has no atoms, 
it follows that for each $i\in I$ we have $\|E_h(\{\lambda_i\})v\|^2=0$, 
i.e.,
$E_h(\{\lambda_i\})v=0$.
On the other hand, since $e_ih=he_i=\lambda_ie_i$, we get $e_i\le
E_h(\{\lambda_i\})$,
that is, $e_iE_h(\{\lambda_i\})=E_h(\{\lambda_i\})e_i=e_i$.
Then $e_iv=e_iE_h(\{\lambda_i\})v=0$ for every $i\in I$.
Since $\sum\limits_{i\in I}e_i=\mathbf{1}$, it then follows that $v=0$, a
contradiction.
\quad $\blacksquare$

\begin{example}\label{newex}
{\rm
A concrete situation where Proposition~\ref{new} applies
is provided by Theorem~2.6.2 in \cite{VDN}.
Specifically, let $M$ be the von Neumann algebra generated by the
real parts $s(t):=(l(t)+l(t)^*)/2$ of the left-creation operators
$l(t)$ (for $t\in{\mathcal H}_{\mathbb R}$) on the full Fock space
${\mathcal T}({\mathcal H}_{\mathbb C})$ associated with
the complexification ${\mathcal H}_{\mathbb C}$ of the real Hilbert space
${\mathcal H}_{\mathbb R}$
with $\dim({\mathcal H}_{\mathbb R})>1$.
Then $M$ is a $\mbox{II}_1$ factor with the trace defined by
the vector form at the vacuum vector $v_0$.

Moreover, for arbitrary $t_0\in{\mathcal H}_{\mathbb R}\setminus\{0\}$,
the operator
$s(t_0)$ is self-adjoint and its spectral measure
localized at the vacuum vector $v_0$
is absolutely continuous with respect to the Lebesgue measure, hence it has no atoms.
As a matter of fact, 
the aforementioned localized spectral measure is given by the semicircle law
$$\frac{2}{\pi\|t_0\|^2}\chi_{[-\|t_0\|,\|t_0\|]}(r)\sqrt{\|t_0\|^2-r^2}\mbox{d}r.$$
Thus, for $\varepsilon>0$ arbitrary, $h:=\|t_0\|+\varepsilon+s(t_0)\in M$ 
is a positive invertible operator
whose spectral measure localized at the vacuum vector $v_0$
(that is, the measure $\|E_h(\cdot)v_0\|^2$)
is absolutely continuous with respect to the Lebesgue measure on ${\mathbb 
R}$.
}
\end{example}

\section{Orbits of adjoint actions}
 Let $\mathcal{ B}(\mathcal{ H})$ denote the space of all bounded 
operators
 and 
${\operatorname{U}}_{\mathcal{ B}(\mathcal{ H})} $ the group of all 
unitary operators
on the complex Hilbert space $\mathcal{H}$.

\begin{theorem}\label{theorem4.1}
{\rm(\cite{DF79})}
Let $\mathcal{ H}$ be a complex Hilbert space,
$T\in\mathcal{ B}(\mathcal{ H})$,  ${\operatorname{U}}_{\mathcal{ 
B}(\mathcal{ H})}(T)$
the unitary
orbit through $T$, and
$$\alpha\colon {\operatorname{U}}_{\mathcal{ B}(\mathcal{ H})}\to
{\operatorname{U}}_{\mathcal{B}(\mathcal{ H})}(T) \;
(\subseteq\mathcal{ B}(\mathcal{ H})), \quad
V\mapsto VTV^*$$
the corresponding orbit map.
The following assertions are equivalent:

\begin{itemize}

\item[{\rm(i)}] The map $\alpha$ has local continuous cross
sections if the unitary orbit ${\operatorname{U}}_{\mathcal{ B}(\mathcal{ 
H})}(T)$ is 
endowed
with the relative topology inherited  from
$\mathcal{B}(\mathcal{H})$.

\item[{\rm(ii)}] The unitary orbit ${\operatorname{U}}_{\mathcal{ 
B}(\mathcal{ H})}(T)$
is closed
in $\mathcal{ B}(\mathcal{ H})$.

\item[{\rm(iii)}] The sub-$C^*$-algebra generated by $T$ in
$\mathcal{ B}(\mathcal{ H})$
is finite dimensional.

\item[{\rm(iv)}] There exist operators $A$ and $B$ on certain
finite dimensional spaces
such that $T$ is unitarily equivalent to the Hilbert space operator
defined by the infinite block-diagonal matrix
$$\left(\begin{array}{cccccc}
A& & &      & &0     \\
 &B& &      & &      \\
 & &B&      & &      \\
 & & &\ddots& &      \\
 & & &      &B&      \\
0& & &      & &\ddots
\end{array}\right).$$

\item[{\rm(v)}] The unitary orbit ${\operatorname{U}}_{\mathcal{ 
B}(\mathcal{ H})}(T)$ 
is a smooth submanifold of $\mathcal{ B}(\mathcal{ H})$.
\end{itemize}
\end{theorem}

\noindent\textbf{Proof.\ \ }
See Theorem~1.1 in \cite{DF79} or Theorem~4.1 in \cite{AFHV84} for the 
fact that 
assertions (i)--(iv) are equivalent.
Moreover, these assertions are equivalent to (v) by the results of 
\cite{AS89}; 
see Theorems 1.1~and~1.3 in \cite{AS91}.
\quad $\blacksquare$

Concerning Theorem~\ref{theorem4.1}(i), we note that the existence of 
global cross sections of the orbit map was investigated in \cite{PH84}.
In fact, according to Theorems 4~and~7 in \cite{PH84}, and using the 
notation of 
Theorem~\ref{theorem4.1},
a global continuous cross section of $\alpha$ can be constructed 
if and only if we can choose $A=B$ in Theorem~\ref{theorem4.1}(iv).

\begin{theorem}\label{theorem4.2}
{\rm(\cite{Ap76})}
Let $\mathcal{ H}$ be a complex Hilbert space, $T\in\mathcal{ B}(\mathcal{
H})$,  and define
$$\operatorname{ad} T\colon\mathcal{ B}(\mathcal{ H})\to\mathcal{ 
B}(\mathcal{
H}) \quad by \quad A\mapsto[T,A].$$
Furthermore denote by ${\mathfrak S}_1$ the ideal of trace class
operators on $\mathcal{H}$.  Then the following assertions are
equivalent:

\begin{itemize}
\item[{\rm(i)}] The operator $\operatorname{ad} T$ has closed range in 
$\mathcal{ 
B}(\mathcal{ H})$.

\item[{\rm(ii)}] For every complex polynomial
$p$  the operator $p(T)$ has closed range in $\mathcal{ H}$ and 
there exists a non-zero polynomial $p_0$ such that $p_0(T)=0$.

\item[{\rm(iii)}] The operator $T$ is similar to an operator
that generates a finite-dimensional sub-$C^*$-algebra of $\mathcal{ 
B}(\mathcal{ H})$.

\item[{\rm(iv)}] The operator $\operatorname{ad} T|_{{\mathfrak S}_1}$ has 
closed range
in
${\mathfrak S}_1$.

\end{itemize}
\end{theorem}

\noindent\textbf{Proof.\ \ }
The fact that the assertions (i)--(iii) are equivalent can be found in
\cite{Ap76}. 

The fact that (i) is equivalent to (iv) is also well known and follows 
by an easy duality argument (see e.g., Theorem~3.5(ii) and the proof of
Proposition~3.12 in \cite{FL84}). 
Thus, first recall that
the range of a Banach space operator is closed if and only if
the range of its dual operator is closed.
Since the Banach space dual to ${\mathfrak S}_1$ is $\mathcal{ 
B}(\mathcal{ H})$ and
the operator dual to $\operatorname{ad} T|_{{\mathfrak S}_1}$ is 
$-\operatorname{ad} T'$
(where $T'\in\mathcal{ B}(\mathcal{ H})$ is the operator dual to $T$),
while $T'$ is conjugate-similar to $T^*$,
it then follows that the range of $\operatorname{ad} T|_{{\mathfrak S}_1}$ 
is closed if
and only if
the range of $\operatorname{ad} T^*$ is closed.
Furthermore, the range of $\operatorname{ad} T^*$ is closed if and only if 
the range of 
$\operatorname{ad} T$
is closed, as a consequence of the fact that 
$\text{(i)}\Leftrightarrow\text{(ii)}$.
\quad $\blacksquare$

\begin{remark}\label{remark4.3}
\normalfont
More details on Theorem~\ref{theorem4.2} can be found in Chapter~15
of the book~\cite{AFHV84}.
\end{remark}

\begin{proposition}\label{proposition4.4}
If $T\in{\mathfrak u}_{\mathcal{ B}(\mathcal{ H})}$ and we denote
$${\operatorname{U}}_{\mathcal{ B}(\mathcal{ H}),T}=
\{U\in {\operatorname{U}}_{\mathcal{ B}(\mathcal{ H})}\mid UTU^{-1}=T\},$$
then ${\operatorname{U}}_{\mathcal{ B}(\mathcal{ H}),T}$ is a Lie subgroup 
of 
${\operatorname{U}}_{\mathcal{ B}(\mathcal{ H})}$.
\end{proposition}

\noindent\textbf{Proof.\ \ }
It follows by Theorem~8.12 in \cite{Be02} that 
${\operatorname{U}}_{\mathcal{
B}(\mathcal{ H}),T}$ is a
Banach Lie group  with respect to the topology inherited from
${\operatorname{U}}_{\mathcal{ B}(\mathcal{ H})}$  and with the Lie 
algebra
${\mathfrak u}_{\mathcal{ B}(\mathcal{ H}),T}
:=\operatorname{Ker}(\operatorname{ad}_{{\mathfrak u}_{\mathcal{ 
B}(\mathcal{ H})}}T)$.

So it only remains to check that ${\mathfrak u}_{\mathcal{ B}(\mathcal{ 
H}),T}$ is a split subspace
of ${\mathfrak u}_{\mathcal{ B}(\mathcal{ H})}$, which is well-known.
Just pick an invariant mean $\mathop{\text{LIM}}\limits_{\alpha\to\infty}$
on the (Abelian, hence amenable) group $({\mathbb R},+)$,
and define a continuous linear map
$$E\colon{\mathfrak u}_{\mathcal{ B}(\mathcal{ H})}\to{\mathfrak 
u}_{\mathcal{ B}(\mathcal{ H})}
\text{ with }E^2=E
\text{ and }\operatorname{Ran} E=
\operatorname{Ker}(\operatorname{ad}_{{\mathfrak u}_{\mathcal{ 
B}(\mathcal{H})}}T)$$
in the following way:
for all $S\in{\mathfrak u}_{\mathcal{ B}(\mathcal{ H})}$ and 
$f,g\in\mathcal{ H}$ let
$$(E(S)f\mid g)=\mathop{\text{LIM}}_{\alpha\to\infty}
(S(\exp(\alpha T))f\mid(\exp(\alpha T))g).$$
We recall that $\mathop{\text{LIM}}\limits_{\alpha\to\infty}$ is just a 
suggestive notation
for a positive linear functional
$m\colon\ell^\infty(\mathbb{R},\mathbb{C})\to\mathbb{C}$
satisfying $\|m\|=1$ and
$$(\forall\xi\in\ell^\infty(\mathbb{R},\mathbb{C}))(\forall\alpha\in\mathbb{R})
\quad
m(\xi)=m(\xi_\alpha),$$
where $\ell^\infty(\mathbb{R},\mathbb{C})$ is the commutative 
$C^*$-algebra
of all bounded functions $\xi\colon\mathbb{R}\to\mathbb{C}$,
and $\xi_\alpha(\beta):=\xi(\alpha+\beta)$
whenever $\xi\in\ell^\infty(\mathbb{R},\mathbb{C})$ and 
$\alpha,\beta\in\mathbb{R}$.
The existence of a functional $m$ with the aforementioned properties 
follows
by Theorem~1.2.1 in \cite{Gr69},
and our notation $\mathop{\text{LIM}}\limits_{\alpha\to\infty}$ is then
introduced by
$$(\forall\xi\in\ell^\infty(\mathbb{R},\mathbb{C}))\quad
\mathop{\text{LIM}}_{\alpha\to\infty}\xi(\alpha):=m(\xi).$$

Now the fact that the map
$E\colon{\mathfrak u}_{\mathcal{ B}(\mathcal{ H})}\to
{\mathfrak u}_{\mathcal{ B}(\mathcal{ H})}$
has the properties claimed above follows by Theorem~16(b) in \cite{Ke02}
applied for the unitary representation
$\alpha\mapsto\exp(\alpha T)$
of the Abelian group $(\mathbb{R},+)$.
\quad $\blacksquare$

\section{Symplectic leaves in preduals of operator ideals}

In this section and in the following one, $\mathcal{H}$ stands for a 
separable complex Hilbert space,
and $\operatorname{GL}(\mathcal{H})$ for the set of all invertible bounded 
linear 
operators on $\mathcal{H}$.

\begin{definition}\label{notation5.1}
%\normalfont
Let $\mathcal{ H}$ be a complex Hilbert space and
${\mathfrak F}$ the ideal of all finite-rank operators on $\mathcal{ H}$.
%\begin{itemize}
%\item[{\rm(i)}]
For every two-sided ideal ${\mathfrak I}$ of $\mathcal{ B}(\mathcal{ H})$
we shall use the following notation:
$$\aligned
{\operatorname{U}}_{\mathfrak I}&={\operatorname{U}}_{\mathcal{ 
B}(\mathcal{ 
H})}\cap(\mathbf{1}+{\mathfrak I}) \\
{\mathfrak u}_{\mathfrak I}&={\mathfrak u}_{\mathcal{ B}(\mathcal{ 
H})}\cap{\mathfrak I}.
\endaligned$$

\end{definition}

For later use, we now recall a few facts concerning Banach ideals of 
operators on
the complex Hilbert space ${\mathcal H}$
(see \cite{GK69} and also \cite{DFWW}).

\begin{remark}\label{ideals}
{\rm
\textbf{(i)} By \textit{Banach ideal} we mean a two-sided ideal 
$\mathfrak{I}$ of
${\mathcal B}({\mathcal H})$
equipped with a norm $\|\cdot\|_\mathfrak{I}$ satisfying
$\|T\|\le\|T\|_\mathfrak{I}=\|T^*\|_\mathfrak{I}$ and
$\|ATB\|_\mathfrak{I}\le\|A\|\,\|T\|_\mathfrak{I} \,\|B\|$
whenever $A,B\in{\mathcal B}({\mathcal H})$.
\medskip

\textbf{(ii)} Let $\widehat{c}$ be the vector space of all sequences of 
real
numbers $\{\xi_j\}_{j\ge1}$
such that $\xi_j=0$ for all but finitely many indices.
A \textit{symmetric norming function} is a function 
$\Phi\colon\widehat{c}\to{\mathbb R}$
satisfying the following conditions:
\begin{itemize}
\item[I)] $\Phi(\xi)>0$ whenever $0\ne\xi\in\widehat{c}$,
\item[II)] $\Phi(\alpha\xi)=|\alpha|\Phi(\xi)$
whenever $\alpha\in{\mathbb R}$ and $\xi\in\widehat{c}$,
\item[III)] $\Phi(\xi+\eta)\le\Phi(\xi)+\Phi(\eta)$ whenever
$\xi,\eta\in\widehat{c}$,
\item[IV)] $\Phi((1,0,0,\dots))=1$,
\item[V)] $\Phi(\{\xi_j\}_{j\ge1})=\Phi(\{\xi_{\pi(j)}\}_{j\ge1})$
whenever $\{\xi_j\}_{j\ge1}\in\widehat{c}$ and
$\pi\colon\{1,2,\dots\}\to\{1,2,\dots\}$
is bijective.
\end{itemize}
Any symmetric norming function $\Phi$ gives rise to two Banach ideals
${\mathfrak S}_\Phi$ and ${\mathfrak S}_\Phi^{(0)}$
as follows.
For every bounded sequence of real numbers $\xi=\{\xi_j\}_{j\ge1}$ define
$$\Phi(\xi):=\sup_{n\ge1}\Phi(\xi_1,\xi_2,\dots,\xi_n,0,0,\dots)\in[0,\infty].$$
For all $T\in{\mathcal B}({\mathcal H})$ denote
$$\|T\|_\Phi:=\Phi(\{s_j(T)\}_{j\ge1})\in[0,\infty],$$
where
$s_j(T)=\inf\{\|T-F\|\mid F\in{\mathcal B}({\mathcal H}),\,{\rm 
rank}\,F<j\}$ whenever
$j\ge1$.
With this notation we can define
$$\begin{aligned}
{\mathfrak S}_\Phi&=\{T\in{\mathcal B}({\mathcal 
H})\mid\|T\|_\Phi<\infty\},\\
{\mathfrak S}_\Phi^{(0)}&=\overline{{\mathfrak F}}^{\|\cdot\|_\Phi}\quad
(\subseteq{\mathfrak S}_\Phi),
\end{aligned}$$
that is, ${\mathfrak S}_\Phi^{(0)}$ is the $\|\cdot\|_\Phi$-closure of
the finite-rank operators ${\mathfrak F}$
in ${\mathfrak S}_\Phi$.
Then $\|\cdot\|_\Phi$ is a norm making ${\mathfrak S}_\Phi$ and 
${\mathfrak S}_\Phi^{(0)}$
into Banach ideals
(see \S 4 in Chapter~III in \cite{GK69}).
Actually, every separable Banach ideal equals ${\mathfrak S}_\Phi^{(0)}$
for some symmetric norming function $\Phi$
(see Theorem~6.2 in Chapter~III in \cite{GK69}).
\medskip

\textbf{(iii)}
For every symmetric norming function $\Phi\colon\widehat{c}\to{\mathbb R}$
there exists a unique symmetric norming function
$\Phi^*\colon\widehat{c}\to{\mathbb R}$
such that
$$\Phi^*(\eta)=\sup\left\{
\frac{1}{\Phi(\xi)}\sum\limits_{j=1}^\infty\xi_j\eta_j
\;\Big|\;\xi=\{\xi_j\}_{j\ge1}\in\widehat{c}\ {\rm and}\ 
\xi_1\ge\xi_2\ge\cdots\ge0
\right\}$$
whenever $\eta=\{\eta_j\}_{j\ge1}\in\widehat{c}$ and
$\eta_1\ge\eta_2\ge\cdots\ge0$.
The function $\Phi^*$ is said to be \textit{adjoint} to $\Phi$
and we always have $(\Phi^*)^*=\Phi$
(see Theorem~11.1 in Chapter~III in \cite{GK69}).
For instance, if $1\le p,q\le\infty$, $1/p+1/q=1$,
$\Phi_p(\xi)=\|\xi\|_{\ell^p}$ and $\Phi_q(\xi)=\|\xi\|_{\ell^q}$ whenever
$\xi\in\widehat{c}$,
then $(\Phi_p)^*=\Phi_q$.
If $\Phi$ is any symmetric norming function then the topological dual of 
the Banach space
${\mathfrak S}_\Phi^{(0)}$ is isometrically isomorphic to ${\mathfrak 
S}_{\Phi^*}$
by means of  the duality pairing
$${\mathfrak S}_{\Phi^*}\times{\mathfrak S}_\Phi^{(0)}\to{\mathbb C},\quad
(T,S)\mapsto\operatorname{Tr}(TS)$$
(see Theorems 12.2~and~12.4 in Chapter~III in \cite{GK69}).
}
\end{remark}

\begin{lemma}\label{lemma7}
Let $k$ be a positive integer and
$${\mathfrak F}_k:=\{T\in{\mathfrak F}\mid{\rm rank}\,T\le k\}.$$
Then for every symmetric norming function $\Phi$
the norms $\|\cdot\|_\Phi$ and $\|\cdot\|$ define the same topology on
the set ${\mathfrak F}_k$.
\end{lemma}

\noindent\textbf{Proof.\ \ }
We essentially follow the idea of proof of Lemma~2.1 in \cite{Bo03}.
Inequalities~(3.12) in Chapter~III in \cite{GK69} show that
$$\xi_1=\Phi_\infty(\xi)\le\Phi(\xi)\le\Phi_1(\xi)=\sum_{j=1}^\infty\xi_j$$
whenever $\xi=\{\xi_j\}_{j\ge1}\in\widehat{c}$ and 
$\xi_1\ge\xi_2\ge\cdots\ge0$.
Since for each $F\in{\mathfrak F}_{2k}$ we have 
$s_{2k+1}(F)=s_{2k+2}(F)=\cdots=0$,
we get
$$(\forall F\in{\mathfrak F}_{2k})\quad
\|F\|=\|F\|_{\Phi_\infty}\le\|F\|_\Phi\le
\|F\|_{\Phi_1}=\sum_{j=1}^\infty s_j(F)\le 2k\cdot s_1(F)=2k\|F\|.$$
On the other hand, the difference of any two operators in ${\mathfrak 
F}_k$
clearly belongs to ${\mathfrak F}_{2k}$,
so that
$$(\forall F_1,F_2\in{\mathfrak F}_k)\quad 
\|F_1-F_2\|\le\|F_1-F_2\|_\Phi\le 2k\|F_1-F_2\|,$$
and the proof is complete.
\quad $\blacksquare$

\begin{lemma}\label{lemma5.2}
Let $\mathcal{ X}_0$ be a reflexive real Banach space and
$A_0\colon\mathcal{X}_0\to\mathcal{X}_0$ a bounded linear operator such
that
$\sup\limits_{t\in{\mathbb R}}\|\exp(tA_0)\|<\infty$.
Then ${\mathcal X}_0=\operatorname{Ker} 
A_0\oplus\overline{\operatorname{Ran} A_0}$.
\end{lemma}

\noindent\textbf{Proof.\ \ }
First endow the complexified space
$\mathcal{ X}:=\mathcal{ X}_0\oplus\text{i}\mathcal{ X}_0$
with a norm making the conjugation
$$C\colon\mathcal{ X}\to\mathcal{ X}, \quad  x+\text{i}y\mapsto x-\text{i}
y,$$
into an isometry (see e.g., Notation~1.1 in \cite{Be02} for a method to
define such a norm).
Thus for all $x,y\in\mathcal{ X}_0$ we have
\begin{equation}\label{**}
\|x\|\le\frac{\|x+{\rm i}y\|+\|x-{\rm i}y\|}{2}=\|x+{\rm i}y\|.
\end{equation}
Then denote by $A\in\mathcal{ B}(\mathcal{ X})$ the unique complex-linear
operator
whose restriction to $\mathcal{ X}_0$ is $A_0$ and commutes with the
conjugation, that is,
$AC=CA$.

On the other hand, denote $M:=\sup\limits_{t\in{\mathbb
R}}\|\exp(tA_0)\|$.
Then for all $z=x+{\rm i}y\in\mathcal{ X}$ and $t\in{\mathbb R}$
we have
$$\|\exp(tA)z\|
=\|\exp(tA_0)x+{\rm i}\exp(tA_0)y\|
\le\|x\|+\|y\|
\le2M\|z\|,$$
where the last inequality follows by \eqref{**}.
Thus $\sup\limits_{t\in{\mathbb R}}\|\exp(tA)\|\le 2M$.
Then it is well known that the norm defined on $\mathcal{ X}$
by $\|z\|_1:=\sup\limits_{t\in{\mathbb R}}\|\exp(tA)z\|$
is equivalent to $\|\cdot\|$ and has the property that
$\|\exp(tA)\|_1=1$
for all $t\in{\mathbb R}$
(see e.g., Lemma~7 in \S 2 in \cite{BD71}).
Since the Banach space $\mathcal{ X}$ is reflexive,
it then follows by Corollary~4.5 in \cite{Ma78}
that
$\mathcal{ 
X}=\operatorname{Ker}(\text{i}A)\oplus\overline{\operatorname{Ran}(\text{i}A)}$,
that is, $\mathcal{ X}=\operatorname{Ker} 
A\oplus\overline{\operatorname{Ran} A}$.

Now, we have $CA=AC$, $\mathcal{ X}_0=\{z\in\mathcal{ X}\mid C(z)=z\}$
and $A|_{\mathcal{ X}_0}=A_0$, so it is straightforward
to show that
${\mathfrak X}_0=\operatorname{Ker} A_0\oplus\overline{\operatorname{Ran} 
A_0}$.
\quad $\blacksquare$

\begin{proposition}\label{proposition5.3}
Let ${\mathfrak I}$ be a Banach ideal whose underlying Banach space is 
reflexive, 
$T\in{\mathfrak u}_{\mathfrak I}$,
and denote
$${\operatorname{U}}_{{\mathfrak J},T}=
\{U\in {\operatorname{U}}_{\mathfrak J}\mid UTU^{-1}=T\}.$$
Then ${\operatorname{U}}_{{\mathfrak J},T}$ is a Lie subgroup of
${\operatorname{U}}_{\mathfrak J}$.
\end{proposition}

\noindent\textbf{Proof.\ \ }
We first recall from Proposition~10.11 in \cite{Be02} that
${\operatorname{U}}_{\mathfrak J}$ is real Banach Lie group whose
Banach Lie algebra is ${\mathfrak u}_{\mathfrak J}$ and that the inclusion 
map
${\operatorname{U}}_{\mathfrak J}\hookrightarrow
{\operatorname{U}}_{\mathcal{ B}(\mathcal{ H})}$
is a homomorphism of Banach Lie groups
(see also Lemma~\ref{lemma1} below).
Since ${\operatorname{U}}_{{\mathfrak J},T}$ is just the inverse image of 
${\operatorname{U}}_{\mathcal{ B}(\mathcal{ H}),T}$
by the aforementioned inclusion map, it follows from
Proposition~\ref{proposition4.4} and Lemma~IV.11 in \cite{Ne00}
that ${\operatorname{U}}_{{\mathfrak J},T}$ is a Banach Lie group with
respect to the 
topology inherited from
${\operatorname{U}}_{\mathfrak J}$ and whose Lie algebra is
${\mathfrak u}_{{\mathfrak J},T}=
\operatorname{Ker}(\operatorname{ad}_{{\mathfrak u}_{\mathfrak J}}T)$.

It only remains to be shown that ${\mathfrak u}_{{\mathfrak J},T}$ is a 
split subspace of
${\mathfrak u}_{\mathfrak J}$.
But this follows by Lemma~\ref{lemma5.2}, since for all $t\in{\mathbb R}$
and $S\in{\mathfrak u}_{\mathfrak J}$
we have
$\bigl(\exp(\operatorname{ad}_{{\mathfrak u}_{\mathfrak J}}tT)\bigr)S
=\text{e}^{tT}S\text{e}^{-tT}$,
whence $\|\exp(\operatorname{ad}_{{\mathfrak u}_{\mathfrak J}}tT)\|\le 1$.
\quad $\blacksquare$

\begin{corollary}\label{corollary5.4}
Let $({\mathfrak B},{\mathfrak J})$ be a pair of Banach ideals whose 
underlying 
Banach spaces are reflexive and assume that the trace pairing
$${\mathfrak B}\times{\mathfrak J}\to{\mathbb C},\quad 
(T,S)\mapsto\operatorname{Tr}(TS)$$
is well defined and
induces a topological isomorphism of the topological dual ${\mathfrak
B}^*$
onto ${\mathfrak J}$.
Then the characteristic distribution of the real Banach Lie-Poisson space
${\mathfrak u}_{\mathfrak B} =({\mathfrak u}_{\mathfrak J})_*$ is
integrable and
all its maximal integral manifolds are symplectic leaves.
\end{corollary}

\noindent\textbf{Proof.\ \ }
The proof is similar to that of Corollaries
\ref{corollary3.8}~and~\ref{corollary3.9},
using Proposition~\ref{proposition5.3} instead of
Proposition~\ref{proposition3.7}.
\quad $\blacksquare$

\begin{example}\label{example5.5}
\normalfont
An obvious example of a pair of Banach ideals $({\mathfrak B},{\mathfrak
J})$
to which Corollary~\ref{corollary5.4} applies is a pair of Schatten ideals
$({\mathfrak S}_p,{\mathfrak S}_q)$
with $p,q\in(1,\infty)$ and $1/p+1/q=1$.
More sophisticated pairs of Banach ideals in duality arise
in the duality theory of operator ideals;
see Remark~\ref{ideals}(iii).
\end{example}

We now consider the problem of constructing invariant complex structures
compatible with the symplectic structures on certain of the leaves
in Corollary~\ref{corollary5.4}.
This problem can be treated by the techniques used in the proof of 
Theorem~VII.6 in \cite{Ne00}.

\begin{proposition}\label{proposition5.6}
Assume that the pair of Banach ideals $({\mathfrak B},{\mathfrak J})$
has the properties that the Banach Lie group 
${\operatorname{U}}_{\mathfrak J}$ is
connected and
the trace pairing
$${\mathfrak B}\times{\mathfrak J}\to{\mathbb C},\quad 
(T,S)\mapsto\operatorname{Tr}(TS)$$
is well defined and
induces a topological isomorphism of the topological dual ${\mathfrak
B}^*$
onto ${\mathfrak J}$.
Let $T\in{\mathfrak u}_{\mathfrak B}\cap{\mathfrak F}$ be a given element
and denote
$${\operatorname{U}}_{{\mathfrak J},T}=\{U\in 
{\operatorname{U}}_{\mathfrak J}\mid UTU^{-1}=T\}.$$
Then the homogeneous space ${\operatorname{U}}_{\mathfrak 
J}/{\operatorname{U}}_{{\mathfrak J},T}$
has a ${\operatorname{U}}_{\mathfrak J}$-invariant weakly K\"ahler structure
and this homogeneous space is weakly immersed into ${\mathfrak
u}_{\mathfrak B}$.
\end{proposition}

\noindent\textbf{Proof.\ \ }
$1^\circ$ {\it Preparations}:
Denote $\sigma(T)=\{\lambda_0,\lambda_1,\dots,\lambda_n\}$ with
$\lambda_0=0$,
and $\mathcal{ H}_i={\operatorname{Ker}}(T-\lambda_i \mathbf{1})$ for 
$i=0,\dots,n$.
Since $T^*=-T$, it follows that we have the orthogonal direct sum
$$\mathcal{ H}=\mathcal{ H}_1\oplus\cdots\oplus\mathcal{
H}_n\oplus\mathcal{ H}_0.$$
Moreover, $\dim\mathcal{ H}_i<\infty$ for $i=1,\dots,n$, since
$T\in{\mathfrak F}$.

Henceforth we will think of the operators on $\mathcal{H}$
as operator matrices with respect to the above orthogonal decomposition.
In particular we have
$$T=
\left(
\begin{array}{cccc}
\lambda_1 &      &         &0          \\
                       &\ddots&         &           \\
                       &      &\lambda_n&           \\
                      0&      &         &\text{\bf 0}
\end{array}
\right),$$
which easily implies that
$$\sigma(\operatorname{ad} T|_{\mathfrak J})=\{\lambda_i-\lambda_j\mid 
0\le i,j\le
n\},$$
and that
\begin{equation}\label{***}
{\mathfrak J}=\bigoplus\limits_{\mu\in\sigma(\operatorname{ad} 
T|_{\mathfrak J})}
\operatorname{Ker}(\operatorname{ad} T|_{\mathfrak J}-\mu).
\end{equation}

$2^\circ$ {\it The isotropy group}:
In particular, it follows that $\operatorname{Ker}(\operatorname{ad} 
T|_{\mathfrak J})$ is 
complemented in ${\mathfrak J}$,
hence the Lie algebra
${\mathfrak u}_{{\mathfrak J},T}=
\operatorname{Ker}(\operatorname{ad} T|_{{\mathfrak u}_{\mathfrak J}})$
of the Lie group ${\operatorname{U}}_{{\mathfrak J},T}$ is complemented in 
${\mathfrak
u}_{\mathfrak J}$.
Since ${\operatorname{U}}_{{\mathfrak J},T}$ is a Lie group with the 
topology inherited
from ${\operatorname{U}}_{\mathfrak J}$
(which follows as in the first part of the proof of
Proposition~\ref{proposition5.3}),
we see that ${\operatorname{U}}_{{\mathfrak J},T}$ is in fact a Lie 
subgroup
of ${\operatorname{U}}_{\mathfrak J}$.

$3^\circ$ {\it The complex structure}:
Since $T^*=-T$, it follows that $\sigma(T)\subseteq{\text i}{\mathbb R}$.
Now we can apply Proposition~8.7 in \cite{Be02} 
with $S={\text i}[0,\infty)$,
${\mathfrak z}={\mathbb R}$ and
$\Psi(\gamma)=\gamma (\operatorname{ad} T|_{{\mathfrak u}_{\mathfrak J}})$
for $\gamma\in{\mathbb R}$
to deduce that the subspace
$${\mathfrak p}:=\bigoplus\limits_{\mu\in(-S)\cap\sigma(\operatorname{ad} 
T|_{\mathfrak 
J})}
\operatorname{Ker}(\operatorname{ad} T|_{\mathfrak J}-\mu)$$
of ${\mathfrak J}$ has the properties
\begin{itemize}
\item[(i)] $[{\mathfrak u}_{{\mathfrak J},T},{\mathfrak
p}]\subseteq{\mathfrak p}$,

\item[(ii)] ${\mathfrak p}\cap\overline{\mathfrak p}
={\mathfrak u}_{{\mathfrak J},T}+\text{i}{\mathfrak u}_{{\mathfrak J},T}$,

\item[(iii)] ${\mathfrak p}+\overline{\mathfrak p}={\mathfrak J}$, and

\item[(iv)] ${\mathfrak p}$ is complemented in ${\mathfrak J}$.
\end{itemize}
Actually it is clear from the expression of ${\mathfrak p}$ that we have
\begin{itemize}
\item[(i')] $V{\mathfrak p}V^{-1}\subseteq{\mathfrak p}$ if $V\in
{\operatorname{U}}_{{\mathfrak J},T}$,
\end{itemize}
\noindent hence Theorem~8.4 in \cite{Be02} shows that there exists a
${\operatorname{U}}_{\mathfrak J}$-invariant complex structure
on the homogeneous space ${\operatorname{U}}_{\mathfrak 
J}/{\operatorname{U}}_{{\mathfrak J},T}$.

$4^\circ$ {\it The symplectic structure}:
Now consider the continuous 2-cocycle of ${\mathfrak u}_{\mathfrak J}$
(actually 2-coboundary) defined by
$T\in({\mathfrak u}_{\mathfrak J})_*\subseteq{\mathfrak u}_{\mathfrak
J}^*$:
$$\omega_T\colon{\mathfrak u}_{\mathfrak J}\times{\mathfrak u}_{\mathfrak
J}\to{\mathbb R},
\quad
\omega_T(X,Y)=\operatorname{Tr}(T[X,Y]).$$
This is just the 2-cocycle that gives rise to the
${\operatorname{U}}_{\mathfrak J}$-invariant weakly symplectic structure
of 
${\operatorname{U}}_{\mathfrak J}/{\operatorname{U}}_{{\mathfrak J},T}$
constructed in Theorem~7.3 in \cite{OR03}.

$5^\circ$ {\it K\"ahler compatibility}:
Note that the above expression of ${\mathfrak p}$ (lower triangular block
matrices,
provided we arrange increasingly the eigenvalues of $T$ on
$\text{i}{\mathbb R}$)
immediately shows that we have $\omega_T({\mathfrak p}\times{\mathfrak
p})=\{0\}$,
that is, ${\mathfrak p}$ is actually a complex polarization of ${\mathfrak
u}_{\mathfrak J}$
relative to the continuous 2-cocycle $\omega_T$ (see e.g., Definition~9.10
in \cite{Be02}).
Furthermore, note that
$${\mathfrak u}_{{\mathfrak J},T}=
\operatorname{Ker}(\operatorname{ad} T|_{{\mathfrak u}_{\mathfrak J}})
=\{X\in{\mathfrak u}_{\mathfrak J}\mid
(\forall Y\in{\mathfrak u}_{\mathfrak J}) \quad \omega_T(X,Y)=0\}.$$
Now a standard reasoning (see e.g., page~77 in \cite{Ne00})
shows that the complex and weakly symplectic invariant structures on the
homogeneous space ${\operatorname{U}}_{\mathfrak 
J}/{\operatorname{U}}_{{\mathfrak J},T}$ are 
compatible, thus making it into
a weakly pseudo-K\"ahler manifold.
This manifold is actually K\"ahler since for all $Z\in{\mathfrak p}$
we have $-\text{i}\omega_T(Z,Z^*)\ge0$
just as in the proof of Lemma~VII.4 in \cite{Ne00}.
\quad $\blacksquare$

\begin{remark}\label{correction}
{\rm
In connection with the hypothesis of Proposition~\ref{proposition5.6}
we note that if ${\mathfrak J}={\mathcal B}({\mathcal H})$
then ${\operatorname{U}}_{\mathfrak J}={\operatorname{U}}_{{\mathcal 
B}({\mathcal H})}$
is well known
to be connected.
Also, if ${\mathfrak J}$ is a separable Banach ideal, then
the Banach Lie group ${\operatorname{U}}_{\mathfrak J}$ is connected as an 
easy
consequence of Theorem~(B)
in \cite{Pa65} and Lemma~\ref{lemma1} below.
In fact, Theorem~(B) in \cite{Pa65} implies that, for a separable Banach 
ideal
${\mathfrak J}$, 
the Banach Lie group $\operatorname{GL}_{\mathfrak J}$ has the same 
homotopy groups as
the direct limit group
$\operatorname{GL}(\infty,{\mathbb C})=
\lim\limits_{\longrightarrow}\operatorname{GL}(n,{\mathbb C})$,
with respect to the natural embeddings
$\operatorname{GL}(n,{\mathbb 
C})\hookrightarrow\operatorname{GL}(n+1,{\mathbb C})$,
\[
A\mapsto\left(\begin{array}{cc} A & 0 \cr 0 &
1\end{array}\right).
\]
In particular, $\operatorname{GL}_{\mathfrak J}$
is connected. Then Lemma~\ref{lemma1} below easily implies
 that $\operatorname{U}_{\mathfrak J}$  is
connected.

Thus, in the special case when ${\mathfrak J}$ is separable and 
${\mathfrak B}\subseteq{\mathfrak J}$, the conclusion of the above
Proposition~\ref{proposition5.6} also follows by the results in 
Chapter~10 in
\cite{Be02}.
} 
\end{remark}

\begin{remark}\label{Hilbert-Schmidt}
\normalfont
We mention that 
in the special case when in Proposition \ref{proposition5.6} we have 
${\mathfrak B}={\mathfrak F}={\mathfrak S}_2$ (the Hilbert-Schmidt ideal) 
the homogeneous space 
${\operatorname{U}}_{\mathfrak J}/{\operatorname{U}}_{{\mathfrak J},T}$
is always a strongly K\"ahler manifold;
see~\cite{Ne00} for details. 
\end{remark}

\section{Embedded orbits in operator ideals}

The unitary orbits of finite-rank self-adjoint operators 
are embedded submanifolds of $\mathcal{B}(\mathcal{H})$,  
according to the results of Andruchow and Stojanoff 
\cite{AS89}, \cite{AS91} 
(see Theorem~\ref{theorem4.1} above). 
In this section we prove a more general version of a similar result of Bona 
\cite{Bo00}, \cite{Bo03} saying that,
on unitary orbits of finite-rank operators on Hilbert spaces,
the natural quotient topology coincides with the trace-class topology.
This fact actually follows by Theorem~\ref{theorem4.1} above and an easy 
topological remark
(see Lemma~\ref{lemma0} below),
so that a version of Theorem~\ref{theorem4.1} involving operator ideals
will automatically lead to a generalization of the aforementioned result 
in \cite{Bo00}, \cite{Bo03}.
That generalization will concern smaller unitary orbits consisting in 
operators of the form
$V^*TV$, where $V$ runs through the set of all unitary operators belonging 
to
$\mathbf{1}+{\mathfrak I}$, for a suitable operator ideal ${\mathfrak I}$.
Additionally, we provide conditions ensuring the existence of invariant
K\"ahler structures
on these smaller unitary orbits (Theorem~\ref{final}).

We now prepare to establish a version of
Theorem~\ref{theorem4.1}(${\rm(iii)}\Rightarrow{\rm(i)}$)
in the more general setting of operator ideals
(see Theorem~\ref{theorem2} below).
The key idea consists in showing that
the main steps of the proof of Theorem~2.1 in \cite{DF79}
can be carried out in the present setting.

\begin{lemma}\label{lemma1}
Let ${\mathfrak I}$ be a Banach ideal of ${\mathcal B}({\mathcal H})$.
Then $\operatorname{GL}_{\mathfrak I}
:=\operatorname{GL}({\mathcal H})\cap(\mathbf{1}+{\mathfrak I})$
is a complex Banach Lie group,
${\operatorname{U}}_{\mathfrak I}:={\operatorname{U}}_{{\mathcal 
B}({\mathcal 
H})}\cap(\mathbf{1}+{\mathfrak I})$
is a real Lie subgroup of $\operatorname{GL}_{\mathfrak I}$,
${\mathfrak p}_{\mathfrak I}:=\{A\in{\mathfrak I}\mid A=A^*\}$
is a real Banach space with the norm inherited from ${\mathfrak I}$,
and the map
$$\Phi\colon{{\operatorname{U}}}_{\mathfrak I}\times{\mathfrak 
p}_{\mathfrak 
I}\to{\operatorname{GL}}_{\mathfrak I},\quad
(V,A)\mapsto V\operatorname{e}^A,$$
is a diffeomorphism.
\end{lemma}

\noindent\textbf{Proof.\ \ }
For the Lie group structures of $\operatorname{GL}_{\mathfrak I}$
and ${\operatorname{U}}_{\mathfrak
I}$
see e.g., Proposition~10.11 in \cite{Be02}.
We just recall that the topology of $\operatorname{GL}_{\mathfrak I}$
is defined by the metric $(V_1,V_2)\mapsto\|V_1-V_2\|_{\mathfrak I}$,
where $\|\cdot\|_{\mathfrak I}$ is the norm of ${\mathfrak I}$.

The fact that the polar decomposition induces a diffeomorphism of
${\operatorname{U}}_{\mathfrak I}\times{\mathfrak p}_{\mathfrak I}$
onto $\operatorname{GL}_{\mathfrak I}$
follows just as in the special case ${\mathfrak I}={\mathfrak S}_p$
treated in Proposition~A.4 in \cite{Ne00}.
\quad $\blacksquare$

\begin{lemma}\label{lemma2}
Let $\Phi$ be a symmetric norming function and
${\mathfrak I}={\mathfrak S}_\Phi$.
Also let $f\colon[0,1]\to{\mathbb R}$ be a continuous nondecreasing
function
such that $0\le f(t)\le t$ whenever $t\in[0,1]$.
Then for every sequence $\{A_n\}_{n\ge1}$ in ${\mathfrak I}$
with $0\le A_n\le\mathbf{1}$ for all $n\ge1$ and
$\lim\limits_{n\to\infty}\|A_n\|_\Phi=0$
we have $f(A_n)\in{\mathfrak I}$ for all $n\ge1$ and
$\lim\limits_{n\to\infty}\|f(A_n)\|_\Phi=0$.
\end{lemma}

\noindent\textbf{Proof.\ \ }
We first recall from Remark~\ref{ideals}(ii) that
$$(\forall T\in{\mathfrak I})\quad
\|T\|_\Phi=\Phi(\{s_j(T)\}_{j\ge1}).$$
Then for every positive integer $n$ we have
$$\begin{aligned}
\|f(A_n)\|_\Phi
&=\Phi(\{s_j(f(A_n))\}_{j\ge1})\\
&=\Phi(\{f(s_j(A_n))\}_{j\ge1})\qquad\mbox{(since $f$ is nondecreasing)}\\
&\le\Phi(\{s_j(A_n)\}_{j\ge1})\qquad\quad\;\mbox{(since $0\le f(t)\le t$ 
for $t\in[0,1]$)}\\
&=\|A_n\|_\Phi.
\end{aligned}$$
Thus $\|f(A_n)\|_\Phi<\infty$ for all $n\ge1$
and $\lim\limits_{n\to\infty}\|f(A_n)\|_\Phi=0$.
\quad $\blacksquare$

\begin{theorem}\label{theorem2}
Let $\Phi$ be a symmetric norming function, ${\mathfrak I}={\mathfrak
S}_\Phi$,
$T=T^*\in{\mathfrak F}$, ${{\operatorname{U}}}_\mathfrak{I}(T) : = \{V^*TV 
\mid V \in
{{\operatorname{U}}}_\mathfrak{I}\}$, and
$$\pi\colon{{\operatorname{U}}}_{\mathfrak 
I}\to{{\operatorname{U}}}_{\mathfrak I}(T),
\quad V\mapsto V^*TV.$$
Then there exist an open neighborhood ${\mathcal D}$ of $T\in{\mathcal 
B}({\mathcal H})$
and a map
$$\varphi\colon{\mathcal D}\cap{{\operatorname{U}}}_{\mathfrak I}(T)
\to{\operatorname{U}}_{\mathfrak I}$$
such that
\begin{itemize}
\item[\rm(i)] $\varphi$ is continuous when ${\mathcal
D}\cap{{\operatorname{U}}}_{\mathfrak I}(T)$
is equipped with the topology inherited from ${\mathcal B}({\mathcal H})$
and ${\operatorname{U}}_{\mathfrak I}$
is equipped with its Lie group topology defined by the metric
$(V_1,V_2)\mapsto\|V_1-V_2\|_\Phi$, 
and
\item[\rm(ii)] $\pi\circ\varphi
=\operatorname{id}_{{\mathcal D}\cap{{\operatorname{U}}}_{\mathfrak 
I}(T)}$.
\end{itemize}
\end{theorem}

For the proof of this theorem we need some notations, remarks, and lemmas.

\begin{notation}\label{notation3}
{\rm
We now introduce some notation that will be used until the end of the
proof of
Theorem~\ref{theorem2}.
\begin{itemize}
\item[\rm(i)] We denote $\sigma(T)=\{\lambda_1,\dots,\lambda_p\}$, where
$\lambda_p=0$.
\item[\rm(ii)] For $i=1,\dots,p$,
we denote ${\mathcal K}_i=\operatorname{Ker}(T-\lambda_i \mathbf{1})$,
$E_i$ the orthogonal projection of ${\mathcal H}$ onto ${\mathcal K}_i$,
and $e_i$ is a polynomial in one variable with real coefficients such that
$E_i=e_i(T)$.
\item[\rm(iii)] We pick an open neighborhood ${\mathcal D}$ of
$T\in{\mathcal B}({\mathcal H})$
such that
$$\max_{1\le i\le p}\sup_{R\in{\mathcal D}}\|e_i(R)-e_i(T)\|<1.$$
\end{itemize}
}
\end{notation}

\begin{remark}\label{remark4}
{\rm
Let $V\in{\operatorname{U}}_{\mathfrak I}$ with $R:=V^*TV\in{\mathcal D}$.
For $i=1,\dots,p$ we have
$$%\begin{aligned}
\|V^*E_iV-E_i\|
%&
=\|V^*e_i(T)V-e_i(T)\|
=\|e_i(V^*TV)-e_i(T)\|%\\
%&
=\|e_i(R)-e_i(T)\|<1,
%\end{aligned}
$$
whence
\begin{equation}\label{(1)}
\|(E_iVE_i)^*(E_iVE_i)-E_i\|<1.
\end{equation}
On the other hand,
the inequality $\|V^*E_iV-E_i\|<1$ also implies that
$$\|E_i-VE_iV^*\|\le\|V\|\,\|V^*E_iV-E_i\|\,\|V^*\|<1,$$
whence
\begin{equation}\label{(2)}
\|(E_iVE_i)(E_iVE_i)^*-E_i\|<1.
\end{equation}
Now \eqref{(1)}~and~\eqref{(2)} show that
$E_iVE_i|_{{\mathcal K}_i}\in\operatorname{GL}(\mathcal{K}_i)$
and thus we have a polar decomposition
$$E_iVE_i=X_iQ_i$$
with $Q_i=|E_iVE_i|=((E_iVE_i)^*(E_iVE_i))^{1/2}$,
$\operatorname{Ker} X_i=\operatorname{Ker} Q_i={\mathcal K}_i^\perp$,
and $X_i|_{{\mathcal K}_i}\in{\operatorname{U}}({\mathcal K}_i)$.
We will denote
$$\psi(V)=X_1^*+\cdots+X_p^*\in{\operatorname{U}}({\mathcal H})$$
whenever $V\in{\operatorname{U}}_{\mathfrak I}$ is as above
(that is, $V^*TV\in{\mathcal D}$).
}
\end{remark}

\begin{notation}\label{notation5}
{\rm
With the notation $\psi(\cdot)$ introduced in Remark~\ref{remark4},
we define
$$\varphi\colon{\mathcal D}\cap{{\operatorname{U}}}_{\mathfrak I}(T)
\to{{\operatorname{U}}}({\mathcal H}) \quad \text{by} \quad
\varphi(V^*TV)=\psi(V) V,$$
where $V\in{\operatorname{U}}_{\mathfrak I}$ and $V^*TV\in{\mathcal D}$.
}
\end{notation}

\begin{lemma}\label{lemma6}
We have a well-defined map
$$\varphi\colon{\mathcal D}\cap
{{\operatorname{U}}}_{\mathfrak I}(T)\to{\operatorname{U}}_{\mathfrak I}$$
satisfying
$\pi\circ\varphi=\operatorname{id}_{{\mathcal 
D}\cap{{\operatorname{U}}}_{\mathfrak I}(T)}$.
\end{lemma}

\noindent\textbf{Proof.\ \ }
$1^\circ$
Let $R=V^*TV=W^*TW\in{\mathcal D}$ with 
$V,W\in{\operatorname{U}}_{\mathfrak I}$.
Then $WV^*\in{\rm U}_{\mathfrak I}\cap\{T\}'={\rm U}_{\mathfrak
I}\cap\{E_1,\dots,E_p\}'$,
so that (with the notation of Remark~\ref{remark4}) we have
$$E_iWE_i=WV^* E_iVE_i=(WV^*X_i) Q_i,$$
where $WV^*X_i|_{{\mathcal K}_i}\in{\rm U}({\mathcal K}_i)$
and 
$\operatorname{Ker}\,(WV^*X_i)=\operatorname{Ker}\,X_i=\operatorname{Ker}\,Q_i$.
Thus the above equalities actually give the polar decomposition of
$E_iWE_i$,
whence
$$\psi(W)=\sum_{i=1}^pX_i^*VW^*=\psi(V) VW^*.$$
Consequently $\psi(W)W=\psi(V)V$,
and thus the definition of $\varphi(R)$ is independent on the choice of 
$V\in{\operatorname{U}}_\mathfrak{I}$
with $R=V^*TV$.

$2^\circ$
We now check that $\varphi(V^*TV)\in{\operatorname{U}}_\mathfrak{I}$ if 
$V\in{\operatorname{U}}_\mathfrak{I}$ and
$V^*TV\in{\mathcal D}$.
First note that for all $i,j\in\{1,\dots,p\}$ we have
$X_i,X_j\in\{E_1,\dots,E_p\}'$,
hence $X_i^*X_j=X_i^*E_iX_j=\delta_{ij}E_i$
and similarly $X_iX_j^*=\delta_{ij}E_i$,
where $\delta_{ij}$ is the Kronecker symbol.
This implies that $\psi(V)\psi(V)^*=\psi(V)^*\psi(V)=\mathbf{1}$.
Thus, in order to show that 
$\varphi(V)=\psi(V)V\in{\operatorname{U}}_\mathfrak{I}$,
it remains to check that $\psi(V)\in{\operatorname{U}}_\mathfrak{I}$.

To this end, note that
$$\delta(V):=\sum_{i=1}^p E_iVE_i\in\mathbf{1}+\mathfrak{I},$$
since $\sum\limits_{i=1}^p E_i=\mathbf{1}$ and 
$V\in\mathbf{1}+\mathfrak{I}$.
On the other hand,
as noted in Remark~\ref{remark4},
we have $E_iVE_i|_{{\mathcal K}_i}\in\operatorname{GL}({\mathcal K}_i)$
for $i=1,\dots,p$,
hence $\delta(V)\in\operatorname{GL}({\mathcal H})$ which proves that
$$\delta(V)\in\operatorname{GL}_\mathfrak{I}.$$
Since it is easy to see that the equality
$\delta(V)=\psi(V)^*(Q_1+\cdots+Q_p)$
is just the polar decomposition of $\delta(V)$,
it then follows by Lemma~\ref{lemma1} that 
$\psi(V)^*\in{\operatorname{U}}_\mathfrak{I}$.
Thus $\psi(V)\in{\operatorname{U}}_\mathfrak{I}$, as desired.

$3^\circ$
To finish the proof we have to show that, if 
$V\in{\operatorname{U}}_\mathfrak{I}$ and
$V^*TV\in{\mathcal D}$,
then $\pi(\varphi(V^*TV))=V^*TV$.
However, since $\psi(V)\in\{T\}'$ and $\varphi(V^*TV)=\psi(V) V$, we have
$\varphi(V^*TV)^* T\varphi(V^*TV)=V^*TV$, as required.
\quad $\blacksquare$

\begin{lemma}\label{lemma8}
The map 
$\varphi\colon{\mathcal D}\cap
{\operatorname{U}}_{\mathfrak I}(T)\to{\operatorname{U}}_{\mathfrak I}$
 is continuous when ${\mathcal D}\cap{\operatorname{U}}_{\mathfrak I}(T)$
is equipped with the topology inherited from ${\mathcal B}({\mathcal H})$ 
and ${\operatorname{U}}_{\mathfrak
I}$
is equipped with its Lie group topology.
\end{lemma}

\noindent\textbf{Proof.\ \ }
$1^\circ$ Let $\{V_n\}_{n\ge1}$ be a sequence in
${\operatorname{U}}_\mathfrak{I}$ such that
$\lim\limits_{n\to\infty}\|V_n^*TV_n-T\|=0$.
We will prove that
$\lim\limits_{n\to\infty}\|\varphi(V_n^*TV_n)-\mathbf{1}\|_\Phi=0$.

Clearly we may assume that $V_n^*TV_n\in{\mathcal D}$ for all $n\ge1$.
Denote $W_n=\varphi(V_n^*TV_n)$,
so that $W_n^*T_nW_n=V_n^*TV_n$ for all $n\ge1$.
Thus we also have $\lim\limits_{n\to\infty}\|W_n^*TW_n-T\|=0$.
Lemma~\ref{lemma7} implies
$\lim\limits_{n\to\infty}\|W_n^*TW_n-T\|_\Phi=0$.

For $i,j\in\{1,\dots,p\}$ and $i\ne j$ we have
$E_iTE_i=\lambda_iE_i$ and $E_jTE_j=\lambda_j E_j$,
hence $E_i[T,W_n]E_j=(\lambda_i-\lambda_j)E_iW_nE_j$.
Then
$$\|E_iW_nE_j\|_\Phi
\le\frac{\|TW_n-W_nT\|_\Phi}{|\lambda_i-\lambda_j|}
=\frac{\|W_n^*TW_n-T\|_\Phi}{|\lambda_i-\lambda_j|},$$
and thus $\lim\limits_{n\to\infty}\|E_iW_nE_j\|_\Phi=0$.

Now let $i\in\{1,\dots,p-1\}$, that is, $\lambda_i\ne0$.
Then $\lim\limits_{n\to\infty}\|W_n^*TW_n-T\|=0$
implies $\lim\limits_{n\to\infty}\|e_i(W_n^*TW_n)-e_i(T)\|=0$,
hence $\lim\limits_{n\to\infty}\|W_n^*E_iW_n-E_i\|=0$.
As in Remark~\ref{remark4} we get
$\lim\limits_{n\to\infty}\|(E_iW_nE_i)^*(E_iW_nE_i)-E_i\|=0$,
that is,
$\lim\limits_{n\to\infty}\|E_i-(E_iW_nE_i)^2\|=0$.
(Note that $E_iW_nE_i=E_i\varphi(V_n^*TV_n) E_i\ge0$ according to
Notation~\ref{notation5} and Remark~\ref{remark4}.) 
Since $\sup\{{\rm rank}\,(E_iW_nE_i)\mid n\ge1\}\le{\rm rank}\,E_i<\infty$
(here we use $\lambda_i\ne0$),
we get by Lemma~\ref{lemma7}
that $\lim\limits_{n\to\infty}\|E_i-(E_iW_nE_i)^2\|_\Phi=0$.
Now Lemma~\ref{lemma2} applied for the function $f(t)=1-(1-t)^{1/2}$
shows that $\lim\limits_{n\to\infty}\|E_i-E_iW_nE_i\|_\Phi=0$.

Next denote
$A_n=(\mathbf{1}-E_p)W_n(\mathbf{1}-E_p)$, $B_n=(\mathbf{1}-E_p)W_nE_p$, 
$C_n=E_pW_n(\mathbf{1}-E_p)$ and $D_n=E_pW_nE_p$,
so that
$$W_n=\left(\begin{array}{cc}
            A_n & B_n \\
            C_n & D_n
            \end{array}\right)$$
in the sense that $W_n=A_n+B_n+C_n+D_n$.
What we have already proved is that
$\lim\limits_{n\to\infty}(\|A_n-(\mathbf{1}-E_p)\|_\Phi+\|B_n\|_\Phi+\|C_n\|_\Phi)=0$.
Since $W_n^*W_n=\mathbf{1}$, we get $B_n^*B_n+D_n^*D_n=E_p$,
so that $\lim\limits_{n\to\infty}\|D_n^*D_n-E_p\|_\Phi=0$.
In other words, $\lim\limits_{n\to\infty}\|E_p-(E_pW_nE_p)^2\|_\Phi=0$,
whence $\lim\limits_{n\to\infty}\|E_p-E_pW_nE_p\|_\Phi=0$ as above,
by making use of Lemma~\ref{lemma2}.
Consequently
$\lim\limits_{n\to\infty}\|W_n-\mathbf{1}\|_\Phi=0$,
as desired.

$2^\circ$
We now prove that
$\varphi\colon{\mathcal D}\cap{\operatorname{U}}_{\mathfrak I}(T)
\to{\operatorname{U}}_{\mathfrak I}$
is continuous at all points of ${\mathcal 
D}\cap{\operatorname{U}}_{\mathfrak I}$.
Let $\{V_n\}_{n\ge1}$ be a sequence in ${\operatorname{U}}_\mathfrak{I}$ 
and $V\in{\operatorname{U}}_\mathfrak{I}$
such that
$\lim\limits_{n\to\infty}\|V_n^*TV_n-V^*TV\|=0$.
We have to show that
$\lim\limits_{n\to\infty}\|\varphi(V_n^*TV_n)-\varphi(V^*TV)\|_\Phi=0$.

To this end, first note that
$\lim\limits_{n\to\infty}\|VV_n^*TV_nV^*-T\|=0$,
hence
$\lim\limits_{n\to\infty}\|\varphi(VV_n^*TV_nV^*)-\mathbf{1}\|_\Phi=0$
by step~$1^\circ$ of the proof.
On the other hand, the operator
$W_n:=\varphi(VV_n^*TV_nV^*)$ has the property
$W_n^*TW_n=VV_n^*TV_nV^*$,
hence
$V_n^*TV_n=V^*W_n^*TW_nV$,
and thus
$\varphi(V_n^*TV_n)=\varphi(V^*W_n^*TW_nV)=\psi(W_nV)W_nV$.
We have $\lim\limits_{n\to\infty}\|W_nV-V\|_\Phi=0$,
hence it will suffice to show that
$\lim\limits_{n\to\infty}\|\psi(W_nV)-\psi(V)\|_\Phi=0$.

Thus we have to show that the map
$$\psi\colon\pi^{-1}({\mathcal 
D})(\subseteq{{\operatorname{U}}}_\mathfrak{I})\to{{\operatorname{U}}}_\mathfrak{I}, 
\quad
W\mapsto\psi(W),$$
is continuous with respect to the topology of 
${\operatorname{U}}_\mathfrak{I}$.
To see this, recall from step~$2^\circ$ of the proof of Lemma~\ref{lemma6}
that, if $W\in{\operatorname{U}}_\mathfrak{I}$ and $W^*TW\in{\mathcal D}$,
then $\delta(W)=\psi(W)^*|\delta(W)|$ is the polar decomposition of
$\delta(W)\in\operatorname{GL}_\mathfrak{I}$.
Now Lemma~\ref{lemma1} along with the obvious continuity of the map
$\delta\colon\pi^{-1}({\mathcal D})\to\operatorname{GL}_\mathfrak{I}$
imply that the map $\psi \colon\pi^{-1}({\mathcal D}) \rightarrow 
{{\operatorname{U}}}_\mathfrak{I}$ is continuous.
\quad $\blacksquare$

\medbreak

\noindent\textbf{Proof of Theorem~\ref{theorem2}.\ \ }
Just use Lemmas \ref{lemma6}~and~\ref{lemma8}
(see also Notation~\ref{notation3}).
\quad $\blacksquare$

\begin{lemma}\label{lemma0}

Let $U$, $Q$, $Q_1$ be topological spaces,
$p\colon U\to Q$ and $\iota\colon Q\to Q_1$ continuous mappings,
and $p_1:=\iota\circ p$.
Assume that the following conditions are satisfied:

\begin{itemize}
\item[\rm(i)] The map $\iota$ is injective.
\item[\rm(ii)] For every $x_1\in Q_1$ there exist a neighborhood $W_1$ of
$x_1$
and a continuous map $\sigma_1\colon W_1\to U$ such that
$p_1\circ\sigma_1=\operatorname{id}_{W_1}$.
\end{itemize}
Then $\iota$ is a homeomorphism of $Q$ onto $Q_1$.
\end{lemma}

\noindent\textbf{Proof.\ \ }
We have by (ii) that the map $p_1$ is onto. 
Since $\iota\circ p=p_1$, it then follows that $\iota$ is onto as well.
Thus it only remains to show that $\iota^{-1}\colon Q_1\to Q$ is 
continuous. 

To this end, let $x_1\in Q_1$ arbitrary. 
According to hypothesis~(ii), there is a continuous map
$\sigma_1\colon W_1\to U$ on some neighborhood $W_1$ of $x_1$ such that 
$p_1\circ\sigma_1=\operatorname{id}_{W_1}$, that is, $\iota\circ 
p\circ\sigma_1=\operatorname{id}_{W_1}$.
Then $\iota^{-1}|_{W_1}=p\circ\sigma_1$,
hence $\iota^{-1}$ is continuous on the neighborhood~$W_1$ of $x_1$. 
Since $x_1\in Q_1$ was arbitrary,
it follows that $\iota^{-1}$ is continuous on the whole set $Q_1$. 
\quad $\blacksquare$
\medskip

Concerning part~(i) in the statement of the next theorem,
we note that it involves two (completely unrelated to each other)
symmetric norming functions.
On the topological level, this corresponds to the fact that any two
symmetric norming functions
define the same topology (in fact, the norm topology) on any unitary orbit
of a finite-rank operator,
as a consequence of Lemma~\ref{lemma7}.
We should point out that there exist a large variety of symmetric norming 
functions,
defining various types of operator ideals like Schatten, Lorentz, Orlicz 
and so on
(see \cite{DFWW} for a survey of this subject).
By way of illustrating this remark, we recall that we have already
mentioned in Remark~\ref{ideals}(iii) the functions
$\Phi_p(\cdot)=\|\cdot\|_{\ell^p}$ that define the Schatten ideals.
For other concrete symmetric norming functions,
see Example~\ref{last} below.

\begin{theorem}\label{final}
Let $\Phi$ and $\Psi$ be symmetric norming functions, 
$\mathfrak{I}={\mathfrak 
S}_\Phi$,
$T=T^*\in{\mathfrak F}$ and ${{\operatorname{U}}}_\mathfrak{I}(T) : = 
\{V^*TV \mid V \in 
{{\operatorname{U}}}_\mathfrak{I}\}$.
Then the following assertions hold:
\begin{itemize}
\item[{\rm(i)}] The orbit map
$$\pi\colon{{\operatorname{U}}}_\mathfrak{I}\to{\mathfrak F},\quad 
V\mapsto V^*TV,$$
induces a diffeomorphism of the homogeneous space 
${{\operatorname{U}}}_{\mathfrak{I}}/{{\operatorname{U}}}_{\mathfrak{I},T}$
onto the submanifold ${{\operatorname{U}}}_{\mathfrak{I}}(T)$ of 
${\mathfrak S}_\Psi$.
\item[{\rm(ii)}] If moreover $\Psi^*=\Phi$ and the Banach Lie group
${\operatorname{U}}_{\mathfrak I}$ is connected,
then the orbit ${{\operatorname{U}}}_{\mathfrak{I}}(T)$
is a ${\operatorname{U}}_{\mathfrak{I}}$-homogeneous weakly K\"ahler
manifold.
\end{itemize}
\end{theorem}

\noindent\textbf{Proof.\ \ }
(i) We first use Lemma~\ref{lemma0} with 
$U={{\operatorname{U}}}_\mathfrak{I}$,
$Q={{\operatorname{U}}}_{\mathfrak{I}}/{{\operatorname{U}}}_{\mathfrak{I},T}$,
$Q_1={\operatorname{U}}_{\mathfrak{I}}(T)$, 
$p\colon{{\operatorname{U}}}_{\mathfrak{I}}
\to{{\operatorname{U}}}_{\mathfrak{I}}/{{\operatorname{U}}}_{\mathfrak{I},T}$ 
the
quotient map
and 
$\iota\colon{{\operatorname{U}}}_{\mathfrak{I}}/{{\operatorname{U}}}_{\mathfrak{I},T}
\to{{\operatorname{U}}}_{\mathfrak{I}}(T)$ induced by
the orbit map $\pi$, to deduce that the differentiable map $\iota$ is a 
homeomorphism,
hence a diffeomorphism.
Note that condition~(ii) in Lemma~\ref{lemma0} is satisfied as a
consequence of
Theorem~\ref{theorem2}.
In order to prove that ${{\operatorname{U}}}_{\mathfrak{I}}(T)$ is an 
embedded submanifold of 
${\mathfrak S}_\Psi$,
we now show that the weak immersion 
$\iota\colon{{\operatorname{U}}}_{\mathfrak{I}}/{{\operatorname{U}}}_{\mathfrak{I},T}\to{\mathfrak 
S}_\Psi$
is actually an immersion.
To this end note that the range of its differential at the point 
$p(\mathbf{1})\in{{\operatorname{U}}}_{\mathfrak{I}}/{{\operatorname{U}}}_{\mathfrak{I},T}$
is
$$\{[T,Y]\mid Y\in{\mathfrak u}_\mathfrak{I}\}=\{[T,Y]\mid 
Y=-Y^*\in{\mathfrak F}\}=
\{[T,Y]\mid Y=-Y^*\in{\mathfrak S}_1\},$$
and this is a closed complemented subspace of ${\mathfrak S}_\Psi$,
as an easy consequence of
Theorem~\ref{theorem4.2} and Lemma~\ref{lemma7}.

(ii) Just use Proposition~\ref{proposition5.6}
along with Remark~\ref{ideals}(iii)
(see also the equality~\eqref{***} in step~$1^\circ$
in the proof of Proposition~\ref{proposition5.6}).
\quad $\blacksquare$

\begin{example}\label{last}
\normalfont
Let $\Pi=\{\pi_j\}_{j\ge1}$ be a sequence of real numbers satisfying the
conditions
\begin{itemize}
\item[{\rm(i)}] $1=\pi_1\ge\pi_2\ge\cdots>0$, and
\item[{\rm(ii)}] $\sum\limits_{j=1}^\infty\pi_j=\infty$.
\end{itemize}
\noindent Let $\mathcal{K}(\mathcal{H}) $ denote the ideal of compact
operators on $\mathcal{H}$ and define
$$\begin{aligned}
{\mathfrak S}_\pi&=\{A\in\mathcal{ K}(\mathcal{ H})\mid
\|A\|_\pi:=\sum_{j=1}^\infty\pi_j\operatorname{s}_j(A)<\infty\}, \\
{\mathfrak S}_\Pi &=\Bigl\{A\in\mathcal{ B}(\mathcal{ H})\mid
\|A\|_\Pi:=\sup_{n\ge1}\frac{\operatorname{s}_1(A)+\cdots
+\operatorname{s}_n(A)}{\pi_1+\cdots+\pi_n}<\infty\Bigr\},
\end{aligned}$$
where $(\operatorname{s}_j(A))_{j\ge1}$ denotes, as usual,
the sequence of singular numbers of an operator $A\in\mathcal{
B}(\mathcal{ H})$
(see e.g., Remark~\ref{ideals}(ii)).
In other words,
${\mathfrak S}_\pi={\mathfrak S}_{\Phi_\pi}={\mathfrak
S}_{\Phi_\pi}^{(0)}$
and ${\mathfrak S}_\Pi={\mathfrak S}_{\Phi_\Pi}$,
where the symmetric norming functions
$\Phi_\pi,\Phi_\Pi\colon\widehat{c}\to{\mathbb R}$ are defined by
$$\Phi_\pi(\xi)=\sum_{j=1}^\infty \pi_j\xi_j \qquad
\mbox{ and } \qquad
\Phi_\Pi(\xi)=\sup_{n\ge1}\frac{\xi_1+\cdots+\xi_n}{\pi_1+\cdots+\pi_n}$$
whenever $\xi=\{\xi_j\}_{j\ge1}\in\widehat{c}$ and 
$\xi_1\ge\xi_2\ge\cdots\ge0$.
We note that $(\Phi_\pi)^*=\Phi_\Pi$ by the comments preceding 
Theorem~15.2
in \cite{GK69}.
It then follows by Theorem~15.2 in \cite{GK69} that
$({\mathfrak S}_\pi,{\mathfrak S}_\Pi)$
is a pair of Banach ideals satisfying the hypotheses of
Proposition~\ref{proposition5.6}.

If moreover the sequence $\Pi=\{\pi_j\}_{j\ge1}$ is regular, in the sense 
that it satisfies the condition
\begin{itemize}
\item[{\rm(iii)}]
$\sup\limits_{n\ge1}(\sum\limits_{j=1}^n\pi_j)/(n\pi_n)<\infty$,
\end{itemize}
\noindent then we have the equality
${\mathfrak S}_\Pi =\bigl\{A\in\mathcal{ B}(\mathcal{ H})\mid
\operatorname{s}_n(A)=O(\pi_n)
\text{ as }n\to\infty\bigr\}$
according to Theorem~14.2 in \cite{GK69}.

We note that, just as in the special case of the similar pair
$({\mathfrak S}_1,\mathcal{ B}(\mathcal{ H}))$,
the dual space ${\mathfrak S}_\Pi$ is in general a {\it non}-separable 
Banach space
(see Theorem~14.1 in \cite{GK69} and Remark~\ref{ideals}(iii)).

For the sake of completeness, we note that in the case when the sequence 
$\Pi$ is constant,
that is, 
$\pi_1=\pi_2=\cdots=1$,                                          
we get ${\mathfrak S}_\pi={\mathfrak S}_1$ the trace class, 
and ${\mathfrak S}_\Pi={\mathcal B}({\mathcal H})$.

This is precisely the situation when the above Theorem~\ref{final}(i)
reduces to Theorem~2.5 in \cite{Bo03} (a part of its 
proof appears already in \cite{Bo00}).
That is, to get the latter result, we have to apply
Theorem~\ref{final}(i) for ${\mathfrak I}={\mathcal B}({\mathcal H})$,
i.e., $\Phi(\{\xi_j\}_{j\ge1})=\max\limits_{j\ge1}|\xi_j|$
and $\Psi(\{\xi_j\}_{j\ge1})=\sum\limits_{j=1}^\infty|\xi_j|$.
\end{example}

\addcontentsline{toc}{section}{Acknowledgments}
\noindent\textbf{Acknowledgments.}
We thank P.~Bona for a number of useful remarks 
and also for sending us
the preprint of his paper \cite{Bo03} that was one of
our motivations for this work.
We thank the referee, K.-H.~Neeb, A.~Odzijewicz, M.~Rieffel, and G.~Weiss 
for their comments
that influenced some of our presentation. 
The first author was partially supported by grant CERES 3-28/2003, 
and by the Swiss NSF through
the SCOPES Program during a one month visit at the EPFL;
the excellent working conditions provided by EPFL are gratefully 
acknowledged.
The second author was 
partially supported by the European Commission and the Swiss Federal
Government through funding for the Research Training Network
\emph{Mechanics and Symmetry in Europe} (MASIE) as well as the Swiss
National Science Foundation.

\end{document}